\documentclass[a4paper]{amsart}

\usepackage{amsmath,amsthm,amsbsy, amstext, amscd, amsxtra,amsopn, latexsym,amssymb,graphicx,graphics}

\title{\large{the  partial positivity of  the curvature in  Riemannian  symmetric spaces}}

\def\shorttitle{the  partial positivity of symmetric spaces}

\newtheorem{theorem}{Theorem}[section]{\bf}{\it}
\newtheorem{proposition}{Proposition}[section]{\bf}{\it}
\newtheorem{lemma}{Lemma}[section]{\bf}{\it}
\newtheorem{remark}{Remark}[section]{\it}{\rm}
\newtheorem*{Proof}{Proof}{\it}{\rm}

\newcommand{\bslide}{\begin{slide}}
\newcommand{\eslide}{\end{slide}}
\newcommand{\frameb}{\begin{frame}}
\newcommand{\framee}{\end{frame}}
\newcommand{\bitemi}{\begin{itemize}}
\newcommand{\eitemi}{\end{itemize}}
\newcommand{\bcenter}{\begin{center}}
\newcommand{\ecenter}{\end{center}}
\newcommand{\bblock}{\begin{block}}
\newcommand{\eblock}{\end{block}}
\newcommand{\benum}{\begin{enumerate}}
\newcommand{\eenum}{\end{enumerate}}

\newcommand{\riem}{Riemannian }

\newcommand{\bthmm}{\begin{thm}}
\newcommand{\ethmm}{\end{thm}}
\newcommand{\blemm}{\begin{lem}}
\newcommand{\elemm}{\end{lem}}
\newcommand{\bppm}{\begin{prop}}
\newcommand{\eppm}{\end{prop}}
\newcommand{\bprfm}{\begin{Proof}}
\newcommand{\eprfm}{\end{Proof}}
\newcommand{\bcorm}{\begin{cor}}
\newcommand{\ecorm}{\end{cor}}

\newcommand{\bthm}{\begin{theorem}}
\newcommand{\ethm}{\end{theorem}}
\newcommand{\blem}{\begin{lemma}}
\newcommand{\elem}{\end{lemma}}
\newcommand{\bpp}{\begin{proposition}}
\newcommand{\epp}{\end{proposition}}
\newcommand{\bprf}{\begin{proof}}
\newcommand{\eprf}{\end{proof}}
\newcommand{\brem}{\begin{remark}}
\newcommand{\erem}{\end{remark}}

\newcommand{\hq}{\ \ }
\newcommand{\np}{\newpage}
\newcommand{\q}{\quad}
\newcommand{\qq}{\qquad}

\newcommand{\fr}{\frac}

\newcommand{\la}{\lambda}
\newcommand{\lb}{\lambda}
\newcommand{\tm}{\times}

\newcommand{\vep}{\varepsilon}

\newcommand{\al}{\alpha}

\newcommand{\tta}{\theta}

\newcommand{\beq}{\begin{eqnarray}}
\newcommand{\eeq}{\end{eqnarray}}
\newcommand{\beqs}{\begin{eqnarray*}}
\newcommand{\eeqs}{\end{eqnarray*}}
\newcommand{\bal}{\begin{align*}}
\newcommand{\eal}{\end{align*}}
\newcommand{\bale}{\begin{aligned}}
\newcommand{\eale}{\end{aligned}}
\newcommand{\bequ}{\begin{equation}}
\newcommand{\eequ}{\end{equation}}
\newcommand{\bequs}{\begin{equation*}}
\newcommand{\eequs}{\end{equation*}}

\newcommand{\bc}{\begin{center}}
\newcommand{\ec}{\end{center}}

\newcommand{\bcase}{\begin{cases}}
\newcommand{\ecase}{\end{cases}}
\newcommand{\bmat}{\begin{matrix}}
\newcommand{\emat}{\end{matrix}}
\newcommand{\bbm}{\begin{bmatrix}}
\newcommand{\ebm}{\end{bmatrix}}
\newcommand{\bpm}{\begin{pmatrix}}
\newcommand{\epm}{\end{pmatrix}}
\newcommand{\bvm}{\begin{vmatrix}}
\newcommand{\evm}{\end{vmatrix}}

\newcommand{\mbbr}{\mathbb{R}}
\newcommand{\mbbc}{\mathbb{C}}
\newcommand{\mbbz}{\mathbb{Z}}

\newcommand{\mfa}{\mathfrak{a}}

\newcommand{\mfe}{\mathfrak{e}}
\newcommand{\mff}{\mathfrak{f}}
\newcommand{\mfg}{\mathfrak{g}}
\newcommand{\mfh}{\mathfrak{h}}

\newcommand{\mfk}{\mathfrak{k}}

\newcommand{\mfo}{\mathfrak{o}}
\newcommand{\mfp}{\mathfrak{p}}

\newcommand{\mfs}{\mathfrak{s}}

\newcommand{\mfu}{\mathfrak{u}}

\begin{document}

\maketitle

\begin{center}
\author
{Xusheng Liu}\\
\address
{School of Mathematical Sciences,\\
 Fudan University, Shanghai 
200433, China } \\
 \email{xshliu@fudan.edu.cn }
\end{center}

\begin{abstract}
In this paper, we determine the partial positivity(resp.,
negativity) of the curvature of all irreducible Riemannian
symmetric spaces. From the classifications of abstract root
systems  and maximal subsystems, we can give the calculations for
symmetric spaces both in  classical types and in exceptional
types.
\end{abstract}

\section{introduction}
We recall the definition of $s$-positive curvature from \cite{Wu}.
A \riem manifold $M$ has $s$-positive(resp., s-negative) curvature
if and only if  for each $x\in M$ and for any $(s+1)$ orthonormal
vectors  $\{e_0,e_1, \cdots,
e_s\}$ in  $M_x, \sum_{i=1}^sK(e_0,e_i)>0$(resp., $<0$), where
$K(e_0,e_i)$ denotes the sectional curvature of the plane spanned
by $e_0$ and $e_i$. The 1-positive curvature is equivalent to
positive sectional curvature, and (n-1)-positive curvature is
equivalent to positive Ricci curvature. \par The manifolds which
have $s$-positive(or negative) curvature have some topological
implications as well as their geometric properties. These
manifolds were studied by Wu\cite{Wu},Shen\cite{Shz},
Kenmotsu and Xia\cite{Kx1}\cite{Kx2}.\par

Among them it was shown
in \cite{Shz} that if a proper open manifold $M$ has $s$-positive
curvature then $M$ has the homotopy type of a $CW$ complex with
cells each of dimension $\le s-1$. In particular,
$H_i(M,\mbbz)=0$, for $i\ge s$.
\par

Frankel \cite{Fr1} showed that two compact totally geodesic
submanifolds in an n-dimensional complete \riem manifold N of
positive sectional curvature must intersect if the sum of their
dimension is greater than or equal to n; and he proved \cite{Fr2}
that if V is an r-dimensional compact totally geodesic immersed
submanifold of N with $2r>n$, then the homomorphism of fundamental
group $\pi_1(V)\to \pi_1(N)$ is surjective. \par
These results had been
generalized in the case of manifolds with partially  positive
curvature by Kenmotsu and Xia in \cite{Kx2}. They showed in an
n-dimensional complete connected  \riem manifold with
k-nonnegative curvature, let V and W be two complete immersed
totally geodesic submanifolds of dimensions r and s, let one of V
and W be compact and suppose N has k-positive curvature either at
all points of V or at all points of W, if $r+s\ge n+k-1$ then V
and W must intersect. As the same time they also proved that if V
is a r-dimensional totally geodesic submanifold with $2r\ge
n+k-1$, then the homomorphism of fundamental  group $\pi_1(V)\to
\pi_1(N)$ is surjective. These show that there exist topological
obstructions for the existence of higher dimensional totally
geodesic submanifolds in the \riem manifolds with partially
positive curvature.
\par
Those results show that the notion of s-positivity(or
negativity) is more subtle curvature condition. Naturally,
determining all the s-values for irreducible \riem symmetric spaces
of compact type is an interesting problem. This has been started
by Lee\cite{Lee2}. She gave the s-values for the symmetric
spaces of classical types  by using of the matrix representation.\par
What about the exceptional cases? Concerning more and more
physicists pay attention to the exceptional geometry, we calculate
the partial positivity of the curvature of all irreducible
Riemannian symmetric spaces of compact type. If $M$ has s-positive
curvature, the dual of $M$ is a simply connected irreducible
symmetric space of noncompact type which has s-negative curvature. \par
Our method is different from hers and we can  deal with both
classical types as well as the exceptional types. For the
classical types we recover Lee's result in [Lee2]. Our main result
is the following table.
\par
\vspace{0.3cm}
\centerline{Table 1.1}
\vspace{0.3cm}
\centerline{
\begin{tabular}{|c|c|c|c|c|c|c|}
  \hline
 Type &   compact type & rank & dimension & s\\
 \hline
 $AI$ &   $SU(n)/SO(n)$  & $n-1$  & $\frac 12(n-1)(n+2)$  &  $\frac {n(n-1)}2$\\
 \hline
 $AII$ &  $SU(2n)/Sp(n)$  & $n-1$  & $(n-1)(2n+1)$  &  $(n-1)(2n-3)$ \\
 \hline
 $AIII$ & $SU(p+q))/S(U_p\times U_q))$  & $min(p, q)$  & $2pq$  &  $1+2(p-1)(q-1)$ \\
 \hline
 $BDI$ & $SO(p+q))/SO(p)\times SO(q)$ & $min(p, q)$  & $pq$  &    $1+(p-1)(q-1)$\\
 \hline
 $DIII$ & $SO(2n))/U(n)$ & $[\frac 12n ] $  & $n(n-1)$  &  $1+(n-2)(n-3)$ \\
 \hline
 $CI$ & $Sp(n)/U(n)$ & $n$  & $n(n+1)$  &  $1+n(n-1)$ \\
 \hline
 $CII$ & $Sp(p+q))/Sp(p)\times Sp(q)$ & $min(p, q)$  & $4pq$  &  $1+4(p-1)(q-1)$ \\
 \hline
 $EI$ & $(\mfe_{6(-78)},\mfs\mfp(4))$ & $6$  & $42$  &  $26$ \\
 \hline
 $EII$ & $(\mfe_{6(-78)},\mfs\mfu(6)+\mfs\mfu(2))$ & $4$  & $40$  &  $19$ \\
 \hline
 $EIII$ & $(\mfe_{6(-78)},\mfs\mfo(10)+\mbbr)$ & $2$  & $32$  &  $11$ \\
 \hline
 $EIV$ & $(\mfe_{6(-78)},\mff(4))$ & $2$  & $26$  &  $10$ \\
 \hline
 $EV$ & $(\mfe_{7(-133)},\mfs\mfu(8))$ & $7$  & $70$  &  $43$ \\
 \hline
 $EVI$ & $(\mfe_{7(-133)},\mfs\mfo(12)+\mfs\mfu(2))$ & $4$  & $64$  &  $31$ \\
 \hline
 $EVII$ & $(\mfe_{7(-133)},\mfe_6+\mbbr)$ & $3$  & $54$  &  $27$ \\
 \hline
 $EVIII$ & $(\mfe_{8(-248)},\mfs\mfo(16))$ & $8$  & $128$  &  $71$ \\
 \hline
 $EIX$ & $(\mfe_{8(-248)},\mfe_7+\mfs\mfu(2))$ & $4$  & $112$  &  $55$ \\
 \hline
 $FI$ & $(\mff_{4(-52)},\mfs\mfp(3)+\mfs\mfu(2))$ & $4$  & $28$  &  $13$ \\
 \hline
 $FII$ & $(\mff_{4(-52)},\mfs\mfo(9))$ & $1$  & $16$  &  $1$ \\
 \hline
 $G$ & $(\mfg_{2(-14)},\mfs\mfu(2)+\mfs\mfu(2))$ & $2$  & $8$  &  $3$ \\
 \hline
 \end{tabular}
 }
\begin{remark}
We note some exceptional  cases. \\
(i)if $r=rank(M)=1$, then $s=r=1$.\\
(ii)For AIII, p=q=2, s=4.\\
(iii)For BDI, p=q=2, s=3; p=q=3,s=6.
\end{remark}

\par
{\bf Acknowledgements.}
\thanks
{ I am  grateful to my supervisor, Prof. Y.L.Xin, for his constant
guidance. He brought \cite{Lee2} to my notice,
and suggested me the method of Lie algebra.
He also pointed out some
mistakes in the course  of this work. }

\np
\section{Abstract root system and subsystem}

An abstract root system in a finite dimensional real inner product
space $V$ with inner product $<,>$ is a finite set $\Delta$ (whose
element is called a root) of $V-\{0\}$ such that (i) $\Delta$ spans
$V$. (ii)For $\al\in \Delta$, the root reflection $s_{\alpha}(h)=h-a_{h,\alpha}\alpha,$
carry $\Delta$ to itself,
where $h\in V, a_{h,\alpha}=\frac{2<h,\mathfrak{\alpha}>}{<\alpha,\alpha>}$
. (iii) $a_{\beta,\alpha}$ is an integer whenever $\alpha,\beta\in\Delta$
.

$l=dimV$ is called the $rank$ of $\Delta$. The Weyl group $W(\Delta)$
is the subgroup on orthogonal group of $V$ generated by the reflection
$s_{\alpha}$ for $\alpha\in\Delta$. An abstract root system is
said to be reduced if $\alpha\in\Delta$ implies $2\alpha\notin$
$\Delta$. If $\alpha$ is a root and $\frac{1}{2}\alpha$ is not
a root, we say that $\alpha$ is reduced. An abstract root system
$\Delta$ is said to be reducible if $\Delta$ admits an nontrivial
disjoint decomposition $\Delta=\Delta'\cup\Delta''$ with every member
of $\Delta'$ orthogonal to every member of $\Delta''$. We say $\Delta$
is irreducible if it admits no such nontrivial decomposition.

For $\alpha,\beta\in\Delta$, the $\alpha$ string containing $\beta$
is the set of all members of $\Delta\cup\{0\}$ of the form $\beta+k\alpha,k\in\mathbb{Z}$.
In fact, there are no gaps, $-p\leq k\leq q,$ $p\geq0,$ $q\geq0$
, and $a_{\beta,\alpha}=p-q$. The $\alpha$ string containing $\beta$
contains at most four roots.

We can choose a lexicographic ordering so that $\Delta=\Delta^{+}\cup\Delta^{-}$
as disjoint sum of the set of positive roots and the set of negative
roots. A root $\alpha$ is called simple if $\alpha>0$ and $\alpha$
does not decompose as the sum of two positive roots. A simple root is
necessarily reduced. We can choose $l$ simple roots $\alpha_{1},\dots,\alpha_{l}$
which are linearly independent, every root $\alpha$ has the form
$\alpha=\sum_{i=1}^{l}m_{i}(\alpha)\alpha_{i}$ with all $m_{i}(\alpha)$
nonnegative or nonpositive. We call $\Pi=\{\alpha_{1},\cdots,\alpha_{l}\}$
a simple system or a fundamental system.

For an abstract root system we associate a Dynkin diagram.

\blem
(classifications of root systems){[}Hel,Kna{]}. Up to an isomorphism the
irreducible reduced abstract root system are $A_{n},B_{n},C_{n},D_{n},E_{6},E_{7},E_{8},F_{4}$
and $G_{2}$. For an nonreduced abstract root system, the reduced
roots form an abstract reduced root system $\Delta_{s}$, the roots
$\alpha\in\Delta$ such that $2\alpha\notin\Delta$ form a reduced
abstract root system $\Delta_{l},$ the Weyl group of $\Delta,\Delta_{s,}\Delta_{l}$
coincide. Up to an isomorphism the only irreducible abstract root system
that are not reduced are of the form $(BC)_{n}$.
\elem

A subset $\Delta_{1}$ of $\Delta$ is called a subsystem of $\Delta$
if it is closed, i.e., (i)if $\alpha\in\Delta_{1}$, then $-\alpha\in\Delta_{1}$
. (ii) if $\alpha,\beta\in\Delta_{1,}\alpha+\beta\in\Delta,$ then
$\alpha+\beta\in\Delta_{1}$. A subsystem $\Delta_{1}$ is called
maximal if $\Delta_{1}$ is a proper subset of $\Delta$ and $\Delta_{1}$
is not properly contained in any properly subsystem of $\Delta$.
A subsystem $\Delta_{1}$ is called $l-1$ maximal if it is maximal
and $rank(\Delta_{1})=l-1$, where $l=rank(\Delta)$.

\blem (classifications of l-1 maximal subsystems){[}Wal{]}. (1)Let
$\Delta$ be an irreducible reduced root system, the maximal
(properly) subsystem is l maximal or l-1 maximal. (2)Let
$\Delta_{1}$ be a l-1 maximal subsystem, $\Pi_{1}$ be a
fundamental system of $\Delta_{1}$, then there exists a
fundamental system $\Pi$ of $\Delta$ such that
$\Pi_{1}=\Pi\cap\Delta_{1}$. (3)Let $\mu=\sum\mu_{i}\alpha_{i}$
be the highest root of $\Delta$,
$\Delta_{1}=\{\alpha\in\Delta|\alpha=\sum
m_{i}\alpha_{i},m_{1}=0\}$ then $\Delta_{1}$ is l-1 maximal if and
only if  $\mu_{1}=1$. Every l-1 maximal subsystem can be obtained
in this manner. \elem

\section{real semisimple lie algebras and restricted root system}
We adopt the notation as in \cite{Hel}. Let $M=G/K$ be an
irreducible Riemannian symmetric  space of noncompact type with
the Cartan decomposition $\mfg=\mfk+\mfp$.  We can extend the
Cartan involution $\theta$ to be complex linear on the complex
simple Lie algebra $\mfg^{\mbbc}$: \beq \tta|_{\mfk^{\mbbc}}=id,
\qq \tta|_{\mfp^{\mbbc}}=-id. \eeq Up to isomorphism the real
simple Lie algebra  $\mfg$ is uniquely determined by $\mfg^{\mbbc}$ and
$\tta$. \par Let $\mfh=\mfh_{\mfk}+\mfa$ be a Cartan subalgebra of
$\mfg$ with $\mfa$ be maximal in $\mfp$.
$$dim \mfa=rank (M)=r.$$
Let $\Delta=\Delta(\mfg^\mbbc,\mfh^\mbbc)$ be the corresponding
root system of $\mfg^\mbbc$.

Every root takes real value on
$\mathfrak{h_{0}=i\mathfrak{h}_{\mathfrak{k}}+\mathfrak{a}}.$
Through the Killing form of $\mathfrak{g}^{\mathbb{C}}$ we can
embed the root into $\mathfrak{h}_{0}$ by
$$\alpha(h)=B(\alpha,h), \q \mbox{for }\q h\in\mathfrak{h}^{\mathbb{C}}.$$
The restricted roots are the elements in $\mathfrak{a}$ of the form
$$\{\lambda=\alpha'=\frac{1}{2}(\alpha-\theta\alpha)|a\in\Delta\}.$$
\par
All the restricted roots form an abstract root system $\Sigma$
which can be nonreduced.
For $\lambda\in\Sigma$, we denote $\mathfrak{g}_{\lambda}$ the
corresponding root space with multiplicity
$$m_{\lambda}=dim\mathfrak{g}_{\lambda}.$$
\par
We can choose the ordering in $\Sigma$ which is compatible to the ordering
of $\Delta$.\par

The restricted roots and the restricted root spaces have the following
properties:

(i) $\mathfrak{g_{0}}=\mathfrak{a}\oplus\mathfrak{m},
\mathfrak{m}=Z_{\mathfrak{k}}(\mathfrak{a})=\{Z\in\mfk|[Z,A]=0\q  \mbox{for all } A\in \mfa \}.
$
.

(ii) $\mathfrak{[g}_{\lambda},\mathfrak{g}_{\mu}]\subset\mathfrak{g}_{\lambda+\mu}$
.

(iii) $\theta$ $(\mathfrak{g}_{\lambda})=\mathfrak{g}_{-\lambda}$
, and $\lambda\in\Sigma$ implies $-\lambda\in\Sigma$.

(iv) $\mathfrak{g}=\mathfrak{g}_{0}+\sum_{\lambda\in\Sigma}\mathfrak{g}_{\lambda}$
as direct sum.

For any $X\in\mathfrak{g}$, there exists $H\in\mathfrak{a},X_{0}\in\mathfrak{m},X_{\lambda}\in\mathfrak{g}_{\lambda}$
so that

$$
X=H+X_{0}+\sum_{\lambda\in\Sigma}X_{\lambda}=(X_{0}+\sum_{\lambda\in\Sigma^{+}}(X_{-\lambda}+\theta X_{-\lambda}))+H+\sum_{\lambda\in\Sigma^{+}}(X_{\lambda}-\theta X_{-\lambda}).
$$

We have Iwasawa decomposition $\mathfrak{g}=\mathfrak{k}\oplus\mathfrak{a}\oplus\sum_{\lambda\in\Sigma^{+}}\mathfrak{g}_{\lambda}.$

Let $h\in\mathfrak{a}$, the centralizer of $h$ in $\mathfrak{p}$
is
\beq\label{zph}
Z_{\mathfrak{p}}(h)=\mathfrak{a}+\sum_{\lambda\in\Sigma^{+},\lambda(h)=0} \mathfrak{g}_{\lambda}.
\eeq
\par
Our objective is to get $s=max_{h\in\mathfrak{a}}\{ dimZ_{\mathfrak{p}}(h)\}$.

$h\in\mathfrak{a}$ is called regular if $\lambda(h)\ne0$ for all
$\lambda\in\Sigma$, otherwise singular. Let $C^{+}=\{ h\in\mathfrak{a}|\lambda(h)>0\; for\;\lambda\in\Sigma^{+}\}$
be a restricted Weyl chamber, its closure is the closed Weyl chamber
$\overline{C^{+}}$. Let $W(\Sigma)$ be the Weyl group of $\Sigma$
, we know for any $h\in\mathfrak{a}$, there exists $w\in W(\Sigma$
) such that $w(h)\in\overline{C^{+}}$.

Let $r=dim$ $\mathfrak{a}$ denote the $rank$ of M. If $h$ is regular
then $dimZ_{\mathfrak{p}}(h)=dim(\mathfrak{a})=r.$ Obviously if
$dimZ_{\mathfrak{p}}(h)$ takes maximal value $h$ must be a singular
element. Now let $h$ be singular, set $\Sigma(h)=\{\lambda\in\Sigma|\lambda(h)=B(\lambda,h)=0\}.$
 It is obviously a subsystem of $\Sigma.$ On the other hand, let
$\Sigma_{1}$ is a subsystem of $\Sigma$ with $rank(\Sigma_{1})<r$
, the span of $\Sigma_{1}$ is $V_{1}$, then $dimV_{1}<r$, Let
$h_{0}$ be any nonzero vector in the orthogonal complement of $V_{1}$
in $\mathfrak{a}$, then $\Sigma(h_{0})$ is a subsystem contains
$\Sigma_{1}$. If $rank(\Sigma_{1})=r,$ the orthogonal complement
of $\Sigma_{1}$ in $\mathfrak{a}$ is $0$, there exists no nonzero
vector $h$ in $\mathfrak{a}$ with $\Sigma(h)$ contains $\Sigma_{1}$
. So if $dim\Sigma(h_{0})$ takes the maximum $s$, $\Sigma(h_{0})$
should be a maximal $l-1$ subsystem.

Let $\Pi'=\{\lambda_{1,}\cdots,\lambda_{r}\}$ be the fundamental
system of $\Sigma$. From the lemma of maximal subsystem we see that
$l-1$ maximal subsystem must be of the form of $\Sigma_{k}=\{\lambda\in\Sigma|\lambda=\sum m'_{i}\lambda_{i},\; m'_{k}=0\},$
 where $1\leq k\leq r$.
\par
Since $s$ is the maximum of $dim Z_{\mfp}(h)$ as $h$ varies in $\mfa$, from (\ref{zph}) we have
\beq \label{sval}
s=r+max_{1\leq k\leq r}\sum_{\lambda\in\Sigma_{k}^{+}}m_{\lambda}.
\eeq

We note that  $s=r$ if $r=1$.
The restricted root systems and the Stake diagrams are listed in {[}Hel,
p532-534{]}.

\section{the partial positivity for riemannian symmetric spaces}

Let $M=U/K$ be a simply connected irreducible symmetric spaces of
compact type. Let $\mfu=\mfk+\mfp_*$ be its Cartan decomposition.
We identify $\mfp_*$ with the tangent space of $M$ at $o=eK$. From
the theory of symmetric spaces, up to constants there exists
uniquely $U$-invariant \riem metric on $M$. We fix such one
invariant metric. We know the curvature tensor at $o$ is
$$
R(X,Y)Z=-[[X,Y],Z],\q \mbox{ for all } X,Y,Z \in \mfp_*.
$$
If $X,Y\in \mfp_*$ are orthonormal vectors, then the sectional curvature of the
plane spanned by $X$ and $Y$ is
$$
K(X,Y)=||[X,Y]||^2.
$$
Thus $K(X,Y)=0$ if and only if $[X,Y]=0$. Suppose the $M$ has s-positive curvature.
From definition we have\\
(1)For any s+1 orthonormal vectors $\{X_0,X_1,\cdots,X_s\}$ in $\mfp_*$
$$\sum_{i=1}^sK(X_0,X_i)=\sum_{i=1}^s||[X_0,X_i]||^2>0.$$
(2)There exists s orthonormal vectors $\{X_0,X_1,\cdots,X_{s-1}\}$ in $\mfp_*$
$$\sum_{i=1}^{s-1}K(X_0,X_i)=\sum_{i=1}^{s-1}||[X_0,X_i]||^2=0.$$
This is equivalent to \\
(3)For any $X\in \mfp_*$,  $dim Z_{\mfp_*}(X)<s+1$,
where
$$Z_{\mfp_*}(X)=\{Y\in\mfp_*| [Y,X]=0 \}$$
 is the centralizer of $X$ in $\mfp_*$.  And\\
(4)There exists at least $X_0\in \mfp_*$, $dim Z_{\mfp_*}(X_0)=s$.\\
So we have\\
$$
s=max_{X\in \mfp_*}\{dim Z_{\mfp_*}(X)\}.
$$

\par
We calculate the s in the following program:

Let $\Pi=\{\alpha_{1},\cdots,\alpha_{l}\}$ be the simple system of
the complex simple Lie algebra $\mathfrak{g}^{\mathbb{C}}$.
Every root can be given as the integral linear combination
\begin{eqnarray}
\alpha=m_1(\alpha)\alpha_1+m_2(\alpha)\alpha_2+\cdots+m_l(\alpha)\alpha_l.
\end{eqnarray}

For
$\alpha\in\Delta$, $\alpha'$ is its restriction ( or projection
) to $\mathfrak{a}$. We have
$$\alpha'_{i}=\lambda_{i^{\theta}}, 1\le i^{\tta}\le r, \mbox{ for } 1\leq i\leq l.$$
Then
$$\alpha'=\sum_{i=1}^lm_i\lambda_{i^{\theta}}=\sum_{j=1}^r m'_{j}(\alpha)\lambda_{j}.$$
From this we can get all restricted roots and their
multiplicities. \par
However, in our calculations, we need not know  the multiplicities explicitly.
From (\ref{sval}), for $1\le k\le r$, let
$$\Delta_{k}=\{\alpha\in \Delta | \al'\ne 0,  m'_{k}(\alpha)=0\}.$$
we denote the number of positive roots in $\Delta(k)$ by $s_k'$,
let $s_k=r+s_k'$.
then we have
$$s=max_{1\le k\le r}\{ s_{k}\}.$$

In the following calculation,
we suppose that $r>1$ since $s=r=1$ if $r=1$.
 For four classical complex simple Lie algebras
we imbed the root system into Euclidean spaces.
We adopt the  Dynkin diagrams of complex simple Lie algebras  as in \cite[p476]{Hel}.
The root system $\Delta=\Delta(\mfg^{\mbbc},\mfh^{\mbbc})$ is reduced while
the restricted root system $\Sigma=\Sigma(\mfg,\mfa)$ may be non reduced.
We denote the rank of $\Delta$ and $\Sigma$ by $l$ and $r$ respectively. Note that
our l,r are r,l in the notation of \cite[p532-534]{Hel}. For simple Lie algebra of exception
type, we can get all the roots as listed in \cite{Frv}.
\par
The following formulas are trivial
\beq \label{eqlk}
(l-k)(l-1-k)=l(l-1)+k^2-k(2l-1).\\
(l-k)(l-1-k)-(l-r)(l-1-r)=(r-k)(2l-1-r-k).
\eeq
\par
\blem \label{lemin}
Let $f(t)=t(T-t), T>0, t $ takes integer values $1,2,\cdots,r,\q   r\le T$, if $T-r\ge 1$, then
$$min_{1\le t\le r}{f(t)}=f(1).$$
\elem

\vskip 1cm
The Dynkin diagram of $\mathfrak{a}_l=\mathfrak{sl}(l+1, \mathbb{C})$ is \\
\centerline
{
\begin{picture}(9, 1)
\put(0.125, 0.5){\circle{0.25}}
\put(0.25, 0.5){\line(1, 0){1}}
\put(1.375, 0.5){\circle{0.25}}
\put(1.5, 0.5){\line(1, 0){0.85}}
\put(2.4, 0.4){$\cdots$}
\put(2.85, 0.5){\line(1, 0){0.9}}
\put(3.875, 0.5){\circle{0.25}}
\put(4.0, 0.5){\line(1, 0){1.}}
\put(5.125, 0.5){\circle{0.25}}
\put(0,0){$\alpha_1$}
\put(1.25, 0){$\alpha_2$}
\put(3.75, 0){$\alpha_{l-1}$}
\put(5,0){$\alpha_l$}
\put(5.3,0.4){.}
\end{picture}
}
\par
Let $\varepsilon_j, 1\le j \le l+1$,  be an orthogonal base of $\mathbb{R}^{l+1}$
with $|\varepsilon_j|^2=\frac 1{2(l+1)}$.
The simple root system of $\mathfrak{a}_l$ is
$$\Pi=\{ \alpha_j=\varepsilon_j-\varepsilon_{j+1},  j=1, 2, \cdots,  l \}, \quad |\alpha_j|^2=\frac 1{l+1}. $$
The positive roots are $$\varepsilon_i-\varepsilon_j=\alpha_i+\cdots+\alpha_{j-1},  1\le i<j\le l+1. $$

\par

(1)AI,
M=$SU(n)/SO(n),  \mfg^{\mbbc}=\mathfrak{a}_l=\mathfrak{sl}(n, \mathbb{C}), r=l=n-1$.\par
In this case, the restricted root system $\Sigma=\Delta$.\par
For $1\le k \le r=l, \lb_k=\al_k'$, then $\al'=(\sum m_i\al_i)'=\sum m_i\lb_i.$\par
We have
$$
\Delta_k^+=\{\alpha>0| m_k(\al)=0 \}.
$$\par
Let
$\al=\vep_i-\vep_j=\al_i+\cdots +\al_j, 1\le i<j\le l+1, m_k(\al)=0 $ if and only if
$$i,j\le k \mbox{ or } i,j>k.$$
Then
\beqs
s_k&=&r+\fr{k(k-1)}2+\frac{(l+1-k)(l-k)}{2}\\
&=&-k(l+1-k)+r+\fr {(l+1)l}2.
\eeqs
From Lemma \ref{lemin} we get
$$
s=max_{1\le k\le r}{s_k}=s_1=l+\fr{l(l-1)}2=\fr {l(l+1)}{2}=\fr {n(n-1)}2.
$$

(2)
AII,  $M=SU(2n)/Sp(n), \mathfrak{g}^{\mbbc}=\mathfrak{a}_l,  l=2n-1, r=n-1,l=2r+1$.\par
We have
$$
\al_1'=\al_3'=\cdots=\al_{2r-1}'=\al_{2r+1}'=0; \q
\al_2'=\la_1,\al_4'=\la_2,\cdots, \al_{2r}'=\la_r.
$$
\par
Let $\al=\vep_i -\vep_j, 1\le\i<j\le l+1=2n,
 \al'=0 \mbox{ if and only if } \al=\al_1,\al_3,\cdots,\al_{2n-1}$.\par

For $1\le k\le r=n-1, \Delta_k=\{\al\in \Delta| m_k'(\al)=0\}=\{\al\in \Delta| m_{2k}'(\al)=0\},$
and \\
$$
m_{2k}(\al)=0 \mbox{ if and only if } i,j\le 2k \mbox{ or } j>i>2k.
$$
We see that
\beqs
s_k&=&r+ \fr{2k(2k-1)}{2}+\fr{(l+1-2k)(l-2k)}{2}-(r+1)\\
    &=& -2k(l+1-2k)+\fr{(l+1)l}2 -1.
\eeqs
From Lemma \ref{lemin} we get
$$
s=s_1=\fr{(l-1)(l-2)}{2}=(n-1)(2n-3).
$$
\par

(3)
AIII, we suppose that $p\le q$,
$M=SU(p+q)/S(U_p\times U_q),  \mathfrak{g}=\mathfrak{a}_l,  l=p+q-1,
r=min(p, q)=p, r\le \fr{l+1}2$.\par
We have
$$
\al_i'=\al_{l+1-i}'=\la_i, 1\le i \le r;\q  \al_j'=0, r<j\le l-r.
$$
\par
Let $\al=\vep_i -\vep_j, 1\le i<j\le l+1,
\al'=0$ if and only if $r<i<j\le l+1-r $. \par
$m_k'(\al)=0$ if $m_k(\al)=0$ or $m_{l+1-k}(\al)=0$. \par
Then $m_k'(\al)=0$ implies
$$1\le i<j\le k \mbox{ or } l+1-k<i<j\le l+1 \mbox{ or } k\le i<j \le l+1-k.$$\par
We get
\beqs
s_k&=&r+2\fr{k(k-1)}{2}+\fr{(l+1-2k)(l-2k)}{2}-\fr{(l+1-2r)(l-2r)}{2}\\
&=& k(k-1)+ 2k^2 -k(2l+1)     + \fr{(l+1)l}2      +r-\fr{(l+1-2r)(l-2r)}{2}\\
&=& - k(2l+2-3k)     + \fr{(l+1)l}2      +r-\fr{(l+1-2r)(l-2r)}{2}.\\
\eeqs
\par
If $\fr{2l+2}3 -r\ge 1$, ie., $r\le \fr{2l-1}{3}$ then $max\{s_k\}=s_1.$\par
As $r\le \fr{l+1}2$, we have $r\le \fr{2l-1}{3}$ when  $l\ge 5$.
\par
we have for $l\ge 5$, i.e., $q\ge 3$,
\beqs
s=s_1&=&r+\fr{(l+1-2)(l-2)}{2}-\fr{(l+1-2r)(l-2r)}{2}\\
&=& r+ \fr 12 (2r-2)(2l+1-2r-2)\\
&=& 1+ (r-1)+ (r-1)(2l-1-2r)\\
&=& 1+ 2(r-1)(l-r)\\
&=& 1+2(p-1)(q-1).
\eeqs
\par
If $l=4, r\le \fr 52, r=2$(we suppose $r \ge 2$),
as $1\cdot (10-3)<2 \cdot (10-6), then s_1>s_2, s=s_2$. In this case $p=2,q=3.$ \par
If $l=3, r\le 2, r=2$,
as $1\cdot (8-3)>2\cdot (8-6), then  s_1<s_2, s=s_2=4$. In this case $p=q=2$. \par

 \vskip 1cm
 The Dynkin diagram of $\mathfrak{b}_l=\mathfrak{so}(2l+1, \mathbb{C})$ is \\
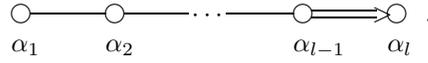
\centerline
{
\begin{picture}(9, 1)
\put(0.125, 0.5){\circle{0.25}}
\put(0.25, 0.5){\line(1, 0){1}}
\put(1.375, 0.5){\circle{0.25}}
\put(1.5, 0.5){\line(1, 0){0.85}}
\put(2.4, 0.4){$\cdots$}
\put(2.85, 0.5){\line(1, 0){0.9}}
\put(3.875, 0.5){\circle{0.25}}
\put(4.0, 0.55){\line(1, 0){0.85}}
\put(4.0, 0.45){\line(1, 0){0.85}}
\put(4.8, 0.4){$>$}
\put(5.125, 0.5){\circle{0.25}}
\put(0, 0){$\alpha_1$}
\put(1.25, 0){$\alpha_2$}
\put(3.75, 0){$\alpha_{l-1}$}
\put(5, 0){$\alpha_l$}
\put(5.5,0.4){.}
\end{picture}
}\par
Let $\varepsilon_j, 1\le j \le l$,  be an orthogonal base of $\mathbb{R}^{l},  |\varepsilon_j|^2=\frac 1{2(2l-1)}$,
the simple root system of $\mathfrak{b}_l$ is
$$\Pi=\{ \alpha_j=\varepsilon_j-\varepsilon_{j+1},  j=1, 2, \cdots,  l-1,  \alpha_l=\varepsilon_l \}. $$
\par The positive roots are $$\varepsilon_i,  \varepsilon_i\pm \varepsilon_j,  1\le i<j \le l.$$
where
\begin{eqnarray*}
&\varepsilon_i=\alpha_i+\cdots+\alpha_l    \\
&\varepsilon_i-\varepsilon_j=\alpha_i+\cdots+\alpha_j \\
&\varepsilon_i+\varepsilon_j=\alpha_i+\cdots+\alpha_{j-1}+2(\alpha_j+\cdots+\alpha_l).
\end{eqnarray*}

(4)BI,
 $M=SO(p+q)/SO(p)\times SO(q), p+q=2l+1$, $\mathfrak{g}=\mathfrak{b}_l, r=p\le l$.\par
 We have
$$
\al_i'=\la_i, 1\le i \le r;\q  \al_j'=0, r<j\le l.
$$\par
If $\al\in \Delta^+,\al'=0$ then $\al$ is one of the following roots,
$$
\vep_i, i>r; \q \vep_i\pm \vep_j, j>i>r.
$$
\par
If $\al>0, m_k'(\al)=0$, then $\al$ is one of the following roots,
$$
\vep_i, i>k; \q
\vep_i - \vep_j, i<j\le k \mbox{ or }  j>i>k; \q
\vep_i + \vep_j, j>i>k.
$$
\par
 Then
 \beqs
 s_k &=&r+ [(l-k)-(l-r)]+ \fr{k(k-1)}{2} +2[\fr{(l-k)(l-1-k)}{2}-\fr{(l-r)(l-1-r)}{2}]\\
 &=& 2r -k + \fr{k(k-1)}{2} + k^2 - k(2l-1) +l(l-1)-(l-r)(l-1-r)\\
 &=& -\fr k2(4l+1-3k)+   2r +l(l-1)-(l-r)(l-1-r).
 \eeqs
 \par
 Since $\fr{4l+1}3-r\ge \fr{4l+1}3-l=\fr{l+1}3\ge 1$, we have
 \beqs
 s=s_1&=&2r-1 + (l-1)(l-2)-(l-r)(l-1-r)\\
 &=& 1+2(r-1)+ (r-1)(2l-1-r-1)\\
 &=& 1+(r-1)(2l-r)\\
 &=& 1+(p-1)(q-1)
 \eeqs
 \par

\vskip 1cm
 The Dynkin diagram of  $\mathfrak{d}_l=\mathfrak{so}(2l, \mathbb{C})$ is \\
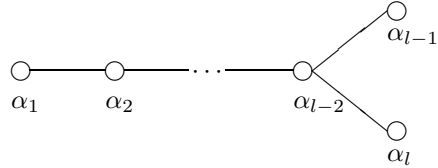
\centerline
{
\begin{picture}(9, 2)
\put(0.125, 0.5){\circle{0.25}}
\put(0.25, 0.5){\line(1, 0){1}}
\put(1.375, 0.5){\circle{0.25}}
\put(1.5, 0.5){\line(1, 0){0.85}}
\put(2.4, 0.4){$\cdots$}
\put(2.85, 0.5){\line(1, 0){0.9}}
\put(3.875, 0.5){\circle{0.25}}
\put(4.0, 0.5){\line(5, 4){1}}
\put(4.0, 0.5){\line(5, -4){1}}
\put(5.125, 1.3){\circle{0.25}}
\put(5.125, -0.3){\circle{0.25}}
\put(5., 0.9){$\alpha_{l-1}$}
\put(5., -0.7){$\alpha_{l}$}
\put(0, 0){$\alpha_1$}
\put(1.25, 0){$\alpha_2$}
\put(3.75, 0){$\alpha_{l-2}$}
\end{picture}
}
\par
\vspace{1cm}
Let $\varepsilon_j, 1\le j \le l$,  be an orthogonal base of $\mathbb{R}^{l}, |\varepsilon_j|^2=\frac 1{4(l-1)}$.
The simple root system is
$$\Pi=\{ \alpha_j=\varepsilon_j-\varepsilon_{j+1},  j=1, 2, \cdots,  l-1,  \alpha_l=\varepsilon_{l-1}+\varepsilon_l \}. $$
The positive roots are
$$\varepsilon_i\pm \varepsilon_j,  i<j. $$
where
\begin{eqnarray*}
&\varepsilon_i-\varepsilon_j&=\alpha_i+\cdots+\alpha_{j-1}\\
&\varepsilon_i+\varepsilon_j&=\alpha_i+\cdots+\alpha_{l-2}+\alpha_{j}+\cdots+\alpha_l\\
&    &=\alpha_i+\cdots+\alpha_{j-1}+2(\alpha_{j}+\cdots+\alpha_{l-2})+\alpha_{l-1}+\alpha_l.
\end{eqnarray*}

(5)DI,
$M=SO(p+q)/SO(p)\times SO(q), p+q=2l$, $\mathfrak{g}=\mathfrak{d}_l, r=p\le l$.\par

There are three cases.\par

(i)$2 \le r\le l-2$,
we have
$$
\al_i'=\la_i, 1\le i \le r;\q  \al_j'=0, r<j\le l.
$$\par
If $\al\in \Delta^+,\al'=0$ then $\al$ is one of the following roots,
$$
\vep_i, i>r; \q \vep_i\pm \vep_j, j>i>r.
$$
\par
If $\al>0, m_k'(\al)=0$, then  $\al$ is one of the following roots,
$$
\vep_i, i>k; \q
\vep_i - \vep_j, i<j\le k \mbox{ or }  j>i>k; \q
\vep_i + \vep_j, j>i>k.
$$
 Then
 \beqs
 s_k &=&r+ \fr{k(k-1)}{2} +2[\fr{(l-k)(l-1-k)}{2}-\fr{(l-r)(l-1-r)}{2}]\\
 &=& \fr{k(k-1)}{2} + k^2 - k(2l-1) +l(l-1)-(l-r)(l-1-r) +r \\
 &=& -\fr k2(4l-1-3k)+   l(l-1)-(l-r)(l-1-r) + r.
 \eeqs
 \par
 Since $\fr{4l-1}3-r\ge \fr{l+5}3\ge 1$, we have
 \beqs
 s=s_1&=&r + (l-1)(l-2)-(l-r)(l-1-r)\\
 &=& 1+(r-1)+ (r-1)(2l-1-r-1)\\
 &=& 1+(r-1)(2l-1-r)\\
 &=& 1+(p-1)(q-1).
 \eeqs\par

 (ii) $r=p=l-1, q=l+1$,
we have
$$
\al_i'=\la_i, 1\le i \le l-2; \al_{l-1}'=\al_l'=\la_{l-1}=\la_r.
$$\par
$\al'>0$ implies $\al=0$.\par
If $\al>0, m_k'(\al)=0$, then  $\al$ is one of the following roots,
$$
\vep_i, i>k; \q
\vep_i - \vep_j, i<j\le k \mbox{ or }  j>i>k; \q
\vep_i + \vep_j, j>i>k.
$$
 Then
 \beqs
 s_k &=&r+ \fr{k(k-1)}{2} +2\fr{(l-k)(l-1-k)}{2}\\
 &=& \fr{k(k-1)}{2} + k^2 - k(2l-1) +l(l-1) +r \\
 &=& -\fr k2(4l-1-3k)+   l(l-1) + r.
 \eeqs
 \par
 Since $\fr{4l-1}3-r=\fr{l+2}3\ge 1$, we have
 \beqs
 s=s_1&=&r + (l-1)(l-2)\\
  &=& (l-1)^2\\
 &=& 1+l(l-2)\\
 &=& 1+(p-1)(q-1).
 \eeqs
 \par

(iii) $r=p=l,q=l$,
we have
$$
\al_i'=\la_i, 1\le i \le l.
$$\par
$\al'>0$ implies $\al=0$.\par
If $\al>0, m_k'(\al)=0$, then  $\al$ is one of the following roots,
$$
\vep_i, i>k; \q
\vep_i - \vep_j, i<j\le k \mbox{ or }  j>i>k; \q
\vep_i + \vep_j, j>i>k.
$$
 Then
 \beqs
 s_k &=&r+ \fr{k(k-1)}{2} +2\fr{(l-k)(l-1-k)}{2}\\
 &=& \fr{k(k-1)}{2} + k^2 - k(2l-1) +l(l-1) +r \\
 &=& -\fr k2(4l-1-3k)+   l(l-1) + r.
 \eeqs
 \par
 Since $\fr{4l-1}3-r=\fr{l-1}3$, we have for $l\ge 4$, i.e., $q\ge 4$,
 \beqs
 s=s_1&=&r + (l-1)(l-2)\\
 &=& 1+(l-1)+ (l-1)(l-2)\\
 &=& 1+(l-1)^2\\
 &=& 1+(p-1)(q-1).
 \eeqs
 \par
When $r\ge 2, l\ge r=2$. For $l=2, r=p=2,q=2, s_1<s_2,s=s_2=3$.\par
For $l=3, r=p=2,q=4, s_1>s_2,s=s_1$;
$r=p=3,q=3, s_3>s_1>s_2,  \\
 s=s_3=6$.

(6)DIII,
$M=SO(2n)/U(n),  \mathfrak{g}^\mbbc=\mathfrak{g}_l,  l=n,r=[\fr l2].$
There are two cases. \par
(i)$l$ is even, we have
$$
\al_{2i}'=\la_i,  1\le i \le r=\fr l2; \al_{2i-1}'=0, 1\le i\le r.
$$\par
$\al>0,\al'>0$ implies $\al=\al_1,\al_3, \cdots, \al_{l-1}$.\par
If $\al>0, m_k'(\al)=0$, then  $\al$ is one of the following roots,
$$
\vep_i - \vep_j, i<j\le 2k \mbox{ or }  j>i>2k; \q
\vep_i + \vep_j, j>i>2k.
$$
 Then
 \beqs
 s_k &=&r+ [\fr{2k(2k-1)}{2} +\fr{(l-2k)(l-1-2k)}{2}-\fr l2]+\fr{(l-2k)(l-1-2k)}2\\
 &=& \fr{k(2k-1)}{2} + (l-2k)(l-1-2k)\\
 &=& -k(4l-1-6k)+   l(l-1).
 \eeqs
 \par
 For  $\fr{4l-1}6- \fr l2= \fr {l-1}6\ge 1$, i.e., $l\ge 7$, we have
 \beqs
 s=s_1&=&1+(l-2)(l-3)=1+(n-2)(n-3).
 \eeqs
 \par
For $l=4,r=2, s_1<s_2, s=s_2=6$.\par
For $l=6,r=3, s_2<s_1<s_3, s=s_3=15$.
\par

(ii)$l$ is odd, we have
$$
\al_{2i}'=\la_i,  1\le i \le \fr {l-3}2; \al_{l-1}'=\al_l'=\la_r, r=\fr{l-1}2;
\al_{2i-1}'=0, 1\le i\le r.
$$\par
$\al>0, \al'>0$ implies $\al=\al_1,\al_3, \cdots, \al_{l-2}$.\par
If $\al>0, m_k'(\al)=0$, then  $\al$ is one of the following roots,
$$
\vep_i - \vep_j, i<j\le 2k \mbox{ or }  j>i>2k; \q
\vep_i + \vep_j, j>i>2k.
$$
 Then
 \beqs
 s_k &=&r+ [\fr{2k(2k-1)}{2} +\fr{(l-2k)(l-1-2k)}{2}-r]+\fr{(l-2k)(l-1-2k)}2\\
 &=& \fr{k(2k-1)}{2} + (l-2k)(l-1-2k)\\
 &=& -k(4l-1-6k)+   l(l-1).
 \eeqs
 \par
 For  $\fr{4l-1}6-\fr {l-1}2=\fr{l+2}6\ge 1$, i.e., $l\ge 4$, we have
 \beqs
 s=s_1&=&1+(l-2)(l-3)=1+(n-2)(n-3).
 \eeqs
 \par
For $l=3,r=1, s=r=1$.
\par

\vskip 1cm
The Dynkin diagram of $\mathfrak{c}_l=\mathfrak{sp}(l, \mathbb{C})$ is \\
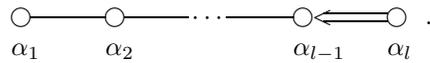
\centerline{
\begin{picture}(9, 1)
\put(0.125, 0.5){\circle{0.25}}
\put(0.25, 0.5){\line(1, 0){1}}
\put(1.375, 0.5){\circle{0.25}}
\put(1.5, 0.5){\line(1, 0){0.85}}
\put(2.4, 0.4){$\cdots$}
\put(2.85, 0.5){\line(1, 0){0.9}}
\put(3.875, 0.5){\circle{0.25}}
\put(4, 0.4){$<$}
\put(4.15, 0.55){\line(1, 0){0.85}}
\put(4.15, 0.45){\line(1, 0){0.85}}
\put(5.125, 0.5){\circle{0.25}}
\put(0, 0){$\alpha_1$}
\put(1.25, 0){$\alpha_2$}
\put(3.75, 0){$\alpha_{l-1}$}
\put(5, 0){$\alpha_l$}
\put(5.5, 0.4){.}
\end{picture}
}
\par
Let $\varepsilon_j, 1\le j \le l$,  be an orthogonal base of
$\mathbb{R}^{l}, |\varepsilon_j|^2=\frac 1{4(l+1)}$,
the simple root system is
$$\Pi=\{ \alpha_j=\varepsilon_j-\varepsilon_{j+1},  j=1, 2, \cdots,  l-1,  \alpha_l=2\varepsilon_l \}. $$
The positive roots are
$$2\varepsilon_i,  \varepsilon_i\pm \varepsilon_j,  i<j.$$
where
\begin{eqnarray*}
&2\varepsilon_i  = & 2(\alpha_i+\cdots+\alpha_{l-1})+\alpha_l\\
&\varepsilon_i-\varepsilon_j  = & \alpha_i+\cdots+\alpha_{j-1}\\
&\varepsilon_i+\varepsilon_j =& \alpha_i+\cdots+\alpha_{j-2}+2(\alpha_{j}+\cdots+\alpha_{l-1})+\alpha_l.
\end{eqnarray*}
\par

(7)CI, $M=Sp(n)/U(n), \mathfrak{g}^\mbbc=\mathfrak{c}_l, l=n$.\par
we have
$$
\al_i'=\la_i,  1\le i \le l.
$$\par

$\al'>0$ implies $\al=0$.\par
If $\al>0, m_k'(\al)=0$, then  $\al$ is one of the following roots,
$$
2\vep_i, i>k; \q
\vep_i - \vep_j, i<j\le 2k \mbox{ or }  j>i>2k; \q
\vep_i + \vep_j, j>i>2k.
$$
 Then
 \beqs
 s_k &=&r+ (l-k) + \fr{k(k-1)}{2} +2\fr{(l-k)(l-1-k)}{2}\\
 &=& 2l-k+\fr{k(k-1)}{2} +(l-k)(l-1-k)\\
 &=& -k(4l+1-3k)+   l(l-1)+2l.
 \eeqs
 \par
 For  $\fr{4l+1}3-l= \fr {l+1}3\ge 1$,i.e., $l\ge 2$, we have
 \beqs
 s=s_1&=&2l-1+(l-1)(l-2)=1+l(l-1)=1+n(n-1).
 \eeqs
 \par

(8)CII,
$M=Sp(p+q)/SP(p)\tm SP(q), \mathfrak{g}^\mbbc=\mathfrak{c}_l, l=p+q, r=p\le [\fr l2]$.\par
We have
$$
\al_{2i}'=\la_i,  1\le i \le r; \al_j'=0, j=1,3,\cdots, 2r-1, \mbox{ or } j>2r.
$$\par

$\al>0,\al'>0$ implies
$\al=\al_1,\al_3, \cdots, \al_{2r-1};
\mbox{ or } \al=\vep_i-\vep_j, j>i>2r$.\par
If $\al>0, m_k'(\al)=0$, then  $\al$ is one of the following roots,
$$
2\vep_i, i>2k; \q
\vep_i - \vep_j, i<j\le 2k \mbox{ or }  j>i>2k; \q
\vep_i + \vep_j, j>i>2k.
$$
Then
 \beqs
 s_k &=&r+ (r-k) + [\fr{2k(2k-1)}{2}-k] +2[\fr{(l-2k)(l-1-2k)}{2}-\fr{(l-2r)(l-1-2r)}{2}]\\
 &=& 2r+ k(2k-3)+(l-2k)(l-1-2k)-(l-2r)(l-1-2r)\\
 &=& -k(4l+1-6k)+ 2r+  l(l-1)- (l-2r)(l-1-2r).
 \eeqs
 \par
 For  $\fr{4l+1}6-\fr l2= \fr {l+1}6\ge 1$,i.e., $l\ge 5$, we have
 \beqs
 s=s_1&=&2r-1 +(l-2)(l-3)-(l-2r)(l-1-2r)\\
 &=&1+4(r-1)(l-1-r)\\
 &=&1+4(p-1)(q-1).
 \eeqs
 \par

For $l=4, r=p=2, q=2, s_1<s_2, s=s_2=6$.\par

\vskip 1cm
Now we  consider the cases of exception type.  \par
We draw the Dynkin diagram  of  $\mfe_6$ as follows. \par
\begin{picture}(9, 2.5)
\put(0.125, 0.5){\circle{0.25}}
\put(0.25, 0.5){\line(1, 0){1}}
\put(1.375, 0.5){\circle{0.25}}
\put(1.5, 0.5){\line(1, 0){1.}}
\put(2.625, 0.5){\circle{0.25}}
\put(2.75, 0.5){\line(1, 0){1.}}
\put(3.875, 0.5){\circle{0.25}}
\put(4.0, 0.5){\line(1, 0){1.}}
\put(5.125, 0.5){\circle{0.25}}
\put(0, 0){$\alpha_1$}
\put(1.25, 0){$\alpha_3$}
\put(2.5, 0){$\alpha_4$}
\put(3.75, 0){$\alpha_5$}
\put(5, 0){$\alpha_6$}
\put(2.5, 0.625){ \line(0, 1){1} }
\put(2.625, 1.75){\circle{0.25}}
\put(2.6, 2){$\alpha_2$}
\put(9, 0.5){\Large{$\mathfrak{e}_6.$}}
\end{picture}\par
We list the positive roots of $\mathfrak{e}_6$ as follows:
\begin{eqnarray*}
&&\alpha_1, \quad \alpha_2, \quad \alpha_3, \quad \alpha_4, \quad \alpha_5, \quad \alpha_6\\
&&\alpha_1+\alpha_3, \quad \alpha_3+\alpha_4, \quad \alpha_2+\alpha_4, \quad\alpha_4+\alpha_5, \quad\alpha_5+\alpha_6\\
&& \alpha_1+\alpha_3+\alpha_4, \quad \alpha_3+\alpha_4+\alpha_2, \quad\alpha_3+\alpha_4+\alpha_5\\
    &&\alpha_4+\alpha_5+\alpha_2, \quad\alpha_4+\alpha_5+\alpha_6\\
&& \alpha_1+\alpha_3+\alpha_4+\alpha_2, \quad \alpha_1+\alpha_3+\alpha_4+\alpha_5, \quad\alpha_2+\alpha_3+\alpha_4+\alpha_5\\
  &&\alpha_3+\alpha_4+\alpha_5+\alpha_6, \quad\alpha_2+\alpha_4+\alpha_5+\alpha_6\\
&& \alpha_1+\alpha_2+\alpha_3+\alpha_4+\alpha_5, \quad \alpha_1+\alpha_3+\alpha_4+\alpha_5+\alpha_6\\
    &&\alpha_2+\alpha_3+\alpha_4+\alpha_5+\alpha_6, \quad\alpha_2+\alpha_3+2\alpha_4+\alpha_5\\
&& \alpha_1+\alpha_2+\alpha_3+\alpha_4+\alpha_5+\alpha_6, \quad \alpha_1+\alpha_2+\alpha_3+2\alpha_4+\alpha_5\\
    &&\alpha_2+\alpha_3+2\alpha_4+\alpha_5+\alpha_6\\
&& \alpha_1+\alpha_2+\alpha_3+2\alpha_4+\alpha_5+\alpha_6, \quad \alpha_1+\alpha_2+2\alpha_3+2\alpha_4+\alpha_5\\
    &&\alpha_2+\alpha_3+2\alpha_4+2\alpha_5+\alpha_6\\
&& \alpha_1+\alpha_2+2\alpha_3+2\alpha_4+\alpha_5+\alpha_6, \quad \alpha_1+\alpha_2+\alpha_3+2\alpha_4+2\alpha_5+\alpha_6\\
&& \alpha_1+\alpha_2+2\alpha_3+2\alpha_4+2\alpha_5+\alpha_6\\
&& \alpha_1+\alpha_2+2\alpha_3+3\alpha_4+2\alpha_5+\alpha_6\\
&& \alpha_1+2\alpha_2+2\alpha_3+3\alpha_4+2\alpha_5+\alpha_6.
\end{eqnarray*}

(9)EI,$\mfg^\mbbc=\mfe_l, l=6,r=6$. \par
we have
$$
\al_i'=\la_i,  1\le i \le 6.
$$\par

$\al'>0$ implies $\al=0$.\par
If $\al>0, m_k'(\al)=0$, then  $m_k(\al)=0$.\par

We list $\Delta_k^+$ as follows,\\
$
\bale
\quad \Delta_1^+=\{&\alpha_{3},\quad    \alpha_{4},\quad    \alpha_{5},\quad    \alpha_{6},\quad    \alpha_{2},\\
    &\alpha_{3}+\alpha_{4},\quad    \alpha_{4}+\alpha_{2},\quad    \alpha_{4}+\alpha_{5},\quad    \alpha_{5}+\alpha_{6},\\
    &     \alpha_{3}+\alpha_{4}+\alpha_{2},\quad    \alpha_{3}+\alpha_{4}+\alpha_{5},\quad    \alpha_{4}+\alpha_{5}+\alpha_{2},\quad    \alpha_{4}+\alpha_{5}+\alpha_{6},\\
    &    \alpha_{3}+\alpha_{4}+\alpha_{5}+\alpha_{2},\quad    \alpha_{3}+\alpha_{4}+\alpha_{5}+\alpha_{6},\quad    \alpha_{4}+\alpha_{5}+\alpha_{6}+\alpha_{2},\\
    &  \alpha_{3}+\alpha_{4}+\alpha_{5}+\alpha_{6}+\alpha_{2},\quad    \alpha_{3}+2\alpha_{4}+\alpha_{5}+\alpha_{2},\\      & \alpha_{3}+2\alpha_{4}+\alpha_{5}+\alpha_{6}+\alpha_{2},
      \alpha_{3}+2\alpha_{4}+2\alpha_{5}+\alpha_{6}+\alpha_{2}  \}\\
\eale
$

$
\bale
\Delta_2^+=\{ & \alpha_{1},\quad    \alpha_{3},\quad    \alpha_{4},\quad    \alpha_{5},\quad    \alpha_{6},\\
    &   \alpha_{1}+\alpha_{3},\quad    \alpha_{3}+\alpha_{4},\quad    \alpha_{4}+\alpha_{5},\quad    \alpha_{5}+\alpha_{6},\\
    & \alpha_{1}+\alpha_{3}+\alpha_{4},\quad    \alpha_{3}+\alpha_{4}+\alpha_{5},\quad    \alpha_{4}+\alpha_{5}+\alpha_{6},\\
    &  \alpha_{1}+\alpha_{3}+\alpha_{4}+\alpha_{5},\quad    \alpha_{3}+\alpha_{4}+\alpha_{5}+\alpha_{6},\\
    &  \alpha_{1}+\alpha_{3}+\alpha_{4}+\alpha_{5}+\alpha_{6}   \}\\
\eale
$

$
\bale
\Delta_3^+=\{ &
\alpha_{1},\quad    \alpha_{4},\quad    \alpha_{5},\quad    \alpha_{6},\quad    \alpha_{2},
      \alpha_{4}+\alpha_{2},\quad    \alpha_{4}+\alpha_{5},\quad    \alpha_{5}+\alpha_{6},\\
    & \alpha_{4}+\alpha_{5}+\alpha_{2},\quad    \alpha_{4}+\alpha_{5}+\alpha_{6},\q
     \alpha_{4}+\alpha_{5}+\alpha_{6}+\alpha_{2}    \}
\eale
$

$
\bale
\Delta_4^+=\{ & \alpha_{1},\quad    \alpha_{3},\quad    \alpha_{5},\quad    \alpha_{6},\quad    \alpha_{2},
&  \alpha_{1}+\alpha_{3},\quad    \alpha_{5}+\alpha_{6}    \}
\eale
$

$
\bale
\Delta_5^+=\{&
\alpha_{1},\quad    \alpha_{3},\quad    \alpha_{4},\quad    \alpha_{6},\quad    \alpha_{2},\q
     \alpha_{1}+\alpha_{3},\quad    \alpha_{3}+\alpha_{4},\quad    \alpha_{4}+\alpha_{2},\\
 &    \alpha_{1}+\alpha_{3}+\alpha_{4},\quad    \alpha_{3}+\alpha_{4}+\alpha_{2},\q
     \alpha_{1}+\alpha_{3}+\alpha_{4}+\alpha_{2}    \}
\eale
$

$
\bale
\Delta_6^+=\{ &
\alpha_{1},\quad    \alpha_{3},\quad    \alpha_{4},\quad    \alpha_{5},\quad    \alpha_{2},\\
 &   \alpha_{1}+\alpha_{3},\quad    \alpha_{3}+\alpha_{4},\quad    \alpha_{4}+\alpha_{2},\quad    \alpha_{4}+\alpha_{5},\\
 &  \alpha_{1}+\alpha_{3}+\alpha_{4},\quad    \alpha_{3}+\alpha_{4}+\alpha_{2},\quad    \alpha_{3}+\alpha_{4}+\alpha_{5},\quad    \alpha_{4}+\alpha_{5}+\alpha_{2},\quad    \alpha_{1}+\alpha_{3}+\alpha_{4}+\alpha_{2},\\
 &   \alpha_{1}+\alpha_{3}+\alpha_{4}+\alpha_{5},\quad    \alpha_{3}+\alpha_{4}+\alpha_{5}+\alpha_{2},\\
 &  \alpha_{1}+\alpha_{3}+\alpha_{4}+\alpha_{5}+\alpha_{2},\quad    \alpha_{3}+2\alpha_{4}+\alpha_{5}+\alpha_{2},\\
 &  \alpha_{1}+\alpha_{3}+2\alpha_{4}+\alpha_{5}+\alpha_{2},\quad
    \alpha_{1}+2\alpha_{3}+2\alpha_{4}+\alpha_{5}+\alpha_{2}    \}.
\eale
$\par

We get the s values as follow\\
$$
s_1=26, s_2=21,s_3=17,s_4=13, s_5=17,s_6=26.$$
Then $s=26$.
\par

(10)EII, $\mfg^\mbbc=\mfe_l, l=6,r=4$. \par
we have
$$
\al_1'=\al_6'=\la_1, \al_3'=\al_5'=\la_3, \al_2'=\la_2, \al_4'=\la_4.
$$\par
$\al'=0$ implies $\al=0$.\par
We list $\Delta_k^+$ as follows,\\
$
\bale
\quad \Delta_1^+=\{ &
\alpha_{3},\quad    \alpha_{4},\quad    \alpha_{5},\quad    \alpha_{2},\\
&   \alpha_{3}+\alpha_{4},\quad    \alpha_{4}+\alpha_{2},\quad    \alpha_{4}+\alpha_{5},\\
&   \alpha_{3}+\alpha_{4}+\alpha_{2},\quad    \alpha_{3}+\alpha_{4}+\alpha_{5},\quad            \alpha_{4}+\alpha_{5}+\alpha_{2},\\
  &   \alpha_{3}+\alpha_{4}+\alpha_{5}+\alpha_{2},\quad    \alpha_{3}+2\alpha_{4}+\alpha_{5}+\alpha_{2}  \}
\eale
$

$
\bale
\Delta_2^+=\{ &
\alpha_{1},\quad    \alpha_{3},\quad    \alpha_{4},\quad    \alpha_{5},\quad    \alpha_{6},\\
 & \alpha_{1}+\alpha_{3},\quad    \alpha_{3}+\alpha_{4},\quad    \alpha_{4}+\alpha_{5},\quad    \alpha_{5}+\alpha_{6},\\
 & \alpha_{1}+\alpha_{3}+\alpha_{4},\quad    \alpha_{3}+\alpha_{4}+\alpha_{5},\quad    \alpha_{4}+\alpha_{5}+\alpha_{6},\\
 &  \alpha_{1}+\alpha_{3}+\alpha_{4}+\alpha_{5},\quad    \alpha_{3}+\alpha_{4}+\alpha_{5}+\alpha_{6},\\
 & \alpha_{1}+\alpha_{3}+\alpha_{4}+\alpha_{5}+\alpha_{6}     \}\\
\eale
$

$
\bale
\Delta_3^+=\{ &
\alpha_{1},\quad    \alpha_{4},\quad    \alpha_{6},\quad    \alpha_{2},\quad    \alpha_{4}+\alpha_{2}     \}
\eale
$

$
\bale
\Delta_4^+=\{
\alpha_{1},\quad    \alpha_{3},\quad    \alpha_{5},\quad    \alpha_{6},\quad    \alpha_{2},\q
 \alpha_{1}+\alpha_{3},\quad    \alpha_{5}+\alpha_{6}     \}.
\eale
$

\par
We get $s_k$ as follows,
$$
s_1=16, s_2=19,s_3=9, s_4=11
.$$
Then $s=19$.
\par

(11)EIII,$\mfg^\mbbc=\mfe_l, l=6,r=2$. \par
we have
$$
\al_1'=\al_6'=\la_1, \al_2'=\la_2, \al_3'=\al_4'=\al_5'=0.
$$\par
$\al'=0$ implies $m_1(\al)=m_6(\al)=0$.\par
We list $\Delta_k^+$ as follows,\\
$
\bale
\quad \Delta_1^+=\{ &
\alpha_{2},\quad    \alpha_{4}+\alpha_{2},\quad    \alpha_{3}+\alpha_{4}+\alpha_{2},\quad    \alpha_{4}+\alpha_{5}+\alpha_{2},\\
&  \alpha_{3}+\alpha_{4}+\alpha_{5}+\alpha_{2},\quad    \alpha_{3}+2\alpha_{4}+\alpha_{5}+\alpha_{2}    \}
\eale
$

$
\bale
\Delta_2^+=\{ &
\alpha_{1},\quad    \alpha_{6},\quad    \alpha_{1}+\alpha_{3},\quad    \alpha_{5}+\alpha_{6},\\
&  \alpha_{1}+\alpha_{3}+\alpha_{4},\quad    \alpha_{4}+\alpha_{5}+\alpha_{6},\\
& \alpha_{1}+\alpha_{3}+\alpha_{4}+\alpha_{5},\quad    \alpha_{3}+\alpha_{4}+\alpha_{5}+\alpha_{6},\\
& \alpha_{1}+\alpha_{3}+\alpha_{4}+\alpha_{5}+\alpha_{6}   \}.
\eale
$\par
We get $s_k$ as follows,
$$
s_1=8, s_2=11
.$$
Then $s=s_2=11$.
\par

(12)EIV,$\mfg^\mbbc=\mfe_l, l=6,r=2$. \par
we have
$$
\al_1'=\la_1, \al_6'=\la_2, \al_2'=\al_3'=\al_4'=\al_5'=0.
$$\par
$\al'=0$ implies $m_1(\al)=m_6(\al)=0$.\par
We list $\Delta_k^+$ as follows,\\
$
\bale
\quad \Delta_1^+=\{ &
\alpha_{6},\quad    \alpha_{5}+\alpha_{6},\quad    \alpha_{4}+\alpha_{5}+\alpha_{6},\\
&  \alpha_{3}+\alpha_{4}+\alpha_{5}+\alpha_{6},\quad    \alpha_{4}+\alpha_{5}+\alpha_{6}+\alpha_{2},\\
&   \alpha_{3}+\alpha_{4}+\alpha_{5}+\alpha_{6}+\alpha_{2},\quad    \alpha_{3}+2\alpha_{4}+\alpha_{5}+\alpha_{6}+\alpha_{2},\\
&   \alpha_{3}+2\alpha_{4}+2\alpha_{5}+\alpha_{6}+\alpha_{2} \}
\eale
$

$
\bale
\Delta_2^+=\{ &
\alpha_{1},\quad    \alpha_{1}+\alpha_{3},\quad    \alpha_{1}+\alpha_{3}+\alpha_{4},\\
  & \alpha_{1}+\alpha_{3}+\alpha_{4}+\alpha_{2},\quad    \alpha_{1}+\alpha_{3}+\alpha_{4}+\alpha_{5},\\
  &   \alpha_{1}+\alpha_{3}+\alpha_{4}+\alpha_{5}+\alpha_{2},\quad    \alpha_{1}+\alpha_{3}+2\alpha_{4}+\alpha_{5}+\alpha_{2},\\
  &   \alpha_{1}+2\alpha_{3}+2\alpha_{4}+\alpha_{5}+\alpha_{2}  \}.
\eale
$\par

We get $s_k$ as follows,
$$
s_1=10, s_2=10
.$$
Then $s=10$.
\par

\vskip 1cm
We draw the Dynkin diagram  of  $\mfe_7$ as follows. \par
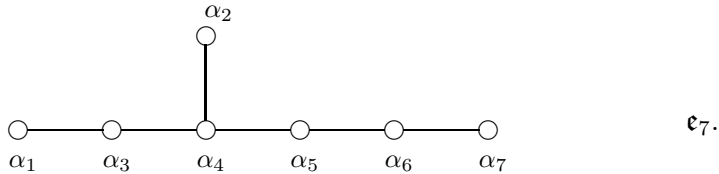
\begin{picture}(9, 2.5)
\put(0.125, 0.5){\circle{0.25}}
\put(0.25, 0.5){\line(1, 0){1}}
\put(1.375, 0.5){\circle{0.25}}
\put(1.5, 0.5){\line(1, 0){1.}}
\put(2.625, 0.5){\circle{0.25}}
\put(2.75, 0.5){\line(1, 0){1.}}
\put(3.875, 0.5){\circle{0.25}}
\put(4.0, 0.5){\line(1, 0){1.}}
\put(5.125, 0.5){\circle{0.25}}
\put(5.25, 0.5){\line(1, 0){1.}}
\put(6.375, 0.5){\circle{0.25}}
\put(0, 0){$\alpha_1$}
\put(1.25, 0){$\alpha_3$}
\put(2.5, 0){$\alpha_4$}
\put(3.75, 0){$\alpha_5$}
\put(5, 0){$\alpha_6$}
\put(6.25, 0){$\alpha_7$}
\put(2.5, 0.625){ \line(0, 1){1} }
\put(2.625, 1.75){\circle{0.25}}
\put(2.6, 2){$\alpha_2$}
\put(9, 0.5){\Large{$\mathfrak{e}_7.$}}
\end{picture}\par

(13)EV,$\mfg^\mbbc=\mfe_l, l=7,r=7$. \par
we have
$$
\al_i'=\la_i, 1\le i\le 7.
$$\par
$\al'=0$ implies $\al=0$.\par
We list $\Delta_k^+$ as follows,\\
$
\bale
\quad \Delta_1^+=\{ &
\alpha_{3},\quad    \alpha_{4},\quad    \alpha_{5},\quad    \alpha_{6},\quad    \alpha_{7},\quad    \alpha_{2},\\
&    \alpha_{3}+\alpha_{4},\quad    \alpha_{4}+\alpha_{5},\quad    \alpha_{4}+\alpha_{2},\quad    \alpha_{5}+\alpha_{6},\quad    \alpha_{6}+\alpha_{7},\\
&   \alpha_{3}+\alpha_{4}+\alpha_{2},\quad    \alpha_{3}+\alpha_{4}+\alpha_{5},\quad    \alpha_{4}+\alpha_{5}+\alpha_{2},\quad    \alpha_{4}+\alpha_{5}+\alpha_{6},\quad    \alpha_{5}+\alpha_{6}+\alpha_{7},\\
 &  \alpha_{3}+\alpha_{4}+\alpha_{5}+\alpha_{2},\quad    \alpha_{3}+\alpha_{4}+\alpha_{5}+\alpha_{6},\quad    \alpha_{4}+\alpha_{5}+\alpha_{6}+\alpha_{2},\\
 &     \alpha_{4}+\alpha_{5}+\alpha_{6}+\alpha_{7},\quad    \alpha_{3}+2\alpha_{4}+\alpha_{5}+\alpha_{2},\\
 &  \alpha_{3}+\alpha_{4}+\alpha_{5}+\alpha_{6}+\alpha_{2},\quad    \alpha_{3}+\alpha_{4}+\alpha_{5}+\alpha_{6}+\alpha_{7},\quad    \alpha_{4}+\alpha_{5}+\alpha_{6}+\alpha_{7}+\alpha_{2},\\
 & \alpha_{3}+2\alpha_{4}+\alpha_{5}+\alpha_{6}+\alpha_{2},\quad    \alpha_{3}+\alpha_{4}+\alpha_{5}+\alpha_{6}+\alpha_{7}+\alpha_{2},\quad    \alpha_{3}+2\alpha_{4}+2\alpha_{5}+\alpha_{6}+\alpha_{2},\\
 &   \alpha_{3}+2\alpha_{4}+\alpha_{5}+\alpha_{6}+\alpha_{7}+\alpha_{2},\quad    \alpha_{3}+2\alpha_{4}+2\alpha_{5}+\alpha_{6}+\alpha_{7}+\alpha_{2},\\
   & \alpha_{3}+2\alpha_{4}+2\alpha_{5}+2\alpha_{6}+\alpha_{7}+\alpha_{2}    \}
\eale
$

$
\bale
\Delta_2^+=\{ &
\alpha_{1},\quad    \alpha_{3},\quad    \alpha_{4},\quad    \alpha_{5},\quad    \alpha_{6},\quad    \alpha_{7},\\
&   \alpha_{1}+\alpha_{3},\quad    \alpha_{3}+\alpha_{4},\quad    \alpha_{4}+\alpha_{5},\quad    \alpha_{5}+\alpha_{6},\quad    \alpha_{6}+\alpha_{7},\\
&    \alpha_{1}+\alpha_{3}+\alpha_{4},\quad    \alpha_{3}+\alpha_{4}+\alpha_{5},\quad    \alpha_{4}+\alpha_{5}+\alpha_{6},\quad    \alpha_{5}+\alpha_{6}+\alpha_{7},\\
&  \alpha_{1}+\alpha_{3}+\alpha_{4}+\alpha_{5},\quad    \alpha_{3}+\alpha_{4}+\alpha_{5}+\alpha_{6},\quad    \alpha_{4}+\alpha_{5}+\alpha_{6}+\alpha_{7},\\
&    \alpha_{1}+\alpha_{3}+\alpha_{4}+\alpha_{5}+\alpha_{6},\quad    \alpha_{3}+\alpha_{4}+\alpha_{5}+\alpha_{6}+\alpha_{7},\quad    \alpha_{1}+\alpha_{3}+\alpha_{4}+\alpha_{5}+\alpha_{6}+\alpha_{7}    \}
\eale
$

$
\bale
\Delta_3^+=\{&
\alpha_{1},\quad    \alpha_{4},\quad    \alpha_{5},\quad    \alpha_{6},\quad    \alpha_{7},\quad    \alpha_{2},\\
&   \alpha_{4}+\alpha_{5},\quad    \alpha_{4}+\alpha_{2},\quad    \alpha_{5}+\alpha_{6},\quad    \alpha_{6}+\alpha_{7},\\
&   \alpha_{4}+\alpha_{5}+\alpha_{2},\quad    \alpha_{4}+\alpha_{5}+\alpha_{6},\quad    \alpha_{5}+\alpha_{6}+\alpha_{7},\\
&    \alpha_{4}+\alpha_{5}+\alpha_{6}+\alpha_{2},\quad    \alpha_{4}+\alpha_{5}+\alpha_{6}+\alpha_{7},\quad    \alpha_{4}+\alpha_{5}+\alpha_{6}+\alpha_{7}+\alpha_{2}    \}
\eale
$

$
\bale
\Delta_4^+=\{ &
\alpha_{1},\quad    \alpha_{3},\quad    \alpha_{5},\quad    \alpha_{6},\quad    \alpha_{7},\quad    \alpha_{2},\\
&  \alpha_{1}+\alpha_{3},\quad    \alpha_{5}+\alpha_{6},\quad    \alpha_{6}+\alpha_{7},\quad    \alpha_{5}+\alpha_{6}+\alpha_{7}     \}
\eale
$

$
\bale
\Delta_5^+=\{&
\alpha_{1},\quad    \alpha_{3},\quad    \alpha_{4},\quad    \alpha_{6},\quad    \alpha_{7},\quad    \alpha_{2},\\
&  \alpha_{1}+\alpha_{3},\quad    \alpha_{3}+\alpha_{4},\quad    \alpha_{4}+\alpha_{2},\quad    \alpha_{6}+\alpha_{7},\\
&   \alpha_{1}+\alpha_{3}+\alpha_{4},\quad    \alpha_{3}+\alpha_{4}+\alpha_{2},\quad    \alpha_{1}+\alpha_{3}+\alpha_{4}+\alpha_{2}    \}
\eale
$

$
\bale
\Delta_6^+=\{ &
\alpha_{1},\quad    \alpha_{3},\quad    \alpha_{4},\quad    \alpha_{5},\quad    \alpha_{7},\quad    \alpha_{2},\\
&    \alpha_{1}+\alpha_{3},\quad    \alpha_{3}+\alpha_{4},\quad    \alpha_{4}+\alpha_{5},\quad    \alpha_{4}+\alpha_{2},\\
&   \alpha_{1}+\alpha_{3}+\alpha_{4},\quad    \alpha_{3}+\alpha_{4}+\alpha_{2},\quad    \alpha_{3}+\alpha_{4}+\alpha_{5},\quad    \alpha_{4}+\alpha_{5}+\alpha_{2},\\
&    \alpha_{1}+\alpha_{3}+\alpha_{4}+\alpha_{2},\quad    \alpha_{1}+\alpha_{3}+\alpha_{4}+\alpha_{5},\quad    \alpha_{3}+\alpha_{4}+\alpha_{5}+\alpha_{2},\quad    \alpha_{3}+2\alpha_{4}+\alpha_{5}+\alpha_{2},\\
&   \alpha_{1}+\alpha_{3}+\alpha_{4}+\alpha_{5}+\alpha_{2},\quad    \alpha_{1}+\alpha_{3}+2\alpha_{4}+\alpha_{5}+\alpha_{2},\quad    \alpha_{1}+2\alpha_{3}+2\alpha_{4}+\alpha_{5}+\alpha_{2}    \}
\eale
$

$
\bale
\Delta_7^+=\{ &
\alpha_{1},\quad    \alpha_{3},\quad    \alpha_{4},\quad    \alpha_{5},\quad    \alpha_{6},\quad    \alpha_{2}, \quad
  \alpha_{1}+\alpha_{3},\quad    \alpha_{3}+\alpha_{4},\quad    \alpha_{4}+\alpha_{5},
    \alpha_{4}+\alpha_{2},\quad    \alpha_{5}+\alpha_{6},\\
&  \alpha_{1}+\alpha_{3}+\alpha_{4},\quad    \alpha_{3}+\alpha_{4}+\alpha_{2},\quad    \alpha_{3}+\alpha_{4}+\alpha_{5},\quad    \alpha_{4}+\alpha_{5}+\alpha_{2}, \quad
  \alpha_{4}+\alpha_{5}+\alpha_{6},\\
&  \alpha_{1}+\alpha_{3}+\alpha_{4}+\alpha_{2},\quad    \alpha_{1}+\alpha_{3}+\alpha_{4}+\alpha_{5},\quad    \alpha_{3}+\alpha_{4}+\alpha_{5}+\alpha_{2},\\
&    \alpha_{3}+\alpha_{4}+\alpha_{5}+\alpha_{6},\quad    \alpha_{4}+\alpha_{5}+\alpha_{6}+\alpha_{2},\quad    \alpha_{3}+2\alpha_{4}+\alpha_{5}+\alpha_{2},\\
&   \alpha_{1}+\alpha_{3}+\alpha_{4}+\alpha_{5}+\alpha_{2},\quad    \alpha_{1}+\alpha_{3}+\alpha_{4}+\alpha_{5}+\alpha_{6},\quad    \alpha_{3}+\alpha_{4}+\alpha_{5}+\alpha_{6}+\alpha_{2},\\
&     \alpha_{1}+\alpha_{3}+2\alpha_{4}+\alpha_{5}+\alpha_{2},\quad    \alpha_{1}+\alpha_{3}+\alpha_{4}+\alpha_{5}+\alpha_{6}+\alpha_{2},\quad    \alpha_{3}+2\alpha_{4}+\alpha_{5}+\alpha_{6}+\alpha_{2},\\
&    \alpha_{1}+2\alpha_{3}+2\alpha_{4}+\alpha_{5}+\alpha_{2},\quad    \alpha_{1}+\alpha_{3}+2\alpha_{4}+\alpha_{5}+\alpha_{6}+\alpha_{2},\quad    \alpha_{3}+2\alpha_{4}+2\alpha_{5}+\alpha_{6}+\alpha_{2},\\
&   \alpha_{1}+2\alpha_{3}+2\alpha_{4}+\alpha_{5}+\alpha_{6}+\alpha_{2}, \quad
    \alpha_{1}+\alpha_{3}+2\alpha_{4}+2\alpha_{5}+\alpha_{6}+\alpha_{2},\\
    &
  \alpha_{1}+2\alpha_{3}+2\alpha_{4}+2\alpha_{5}+\alpha_{6}+\alpha_{2}, \quad
\alpha_{1}+2\alpha_{3}+3\alpha_{4}+2\alpha_{5}+\alpha_{6}+\alpha_{2},\\
&   \alpha_{1}+2\alpha_{3}+3\alpha_{4}+2\alpha_{5}+\alpha_{6}+2\alpha_{2}    \}.
\eale
$
\par
We get $s_k$ as follows,
$$
s_1=37, s_2=28, s_3=23, s_4=17, s_5=20, s_6=28, s_7=43
.$$
Then $s=43$.
\par

(14)EVI,$\mfg^\mbbc=\mfe_l, l=7,r=4$. \par
we have
$$
\al_1'=\la_1, \al_3'=\la_2,\al_4'=\la_3,\al_6'=\la_4,
\al_2'=\al_5'=\al_7'=0.
$$\par
$\al'=0$ implies $m_2(\al)=m_5(\al)=m_7(\al)=0$.\par
We list $\Delta_k^+$ as follows,\\
$
\bale
\quad \Delta_1^+=\{ &
\alpha_{3},\quad    \alpha_{4},\quad    \alpha_{6},\quad    \alpha_{3}+\alpha_{4},\quad    \alpha_{4}+\alpha_{5},\quad    \alpha_{4}+\alpha_{2},\quad    \alpha_{5}+\alpha_{6},\quad    \alpha_{6}+\alpha_{7},\\
&  \alpha_{3}+\alpha_{4}+\alpha_{2},\quad    \alpha_{3}+\alpha_{4}+\alpha_{5},\quad    \alpha_{4}+\alpha_{5}+\alpha_{2},\quad    \alpha_{4}+\alpha_{5}+\alpha_{6},\quad    \alpha_{5}+\alpha_{6}+\alpha_{7},\\
&  \alpha_{3}+\alpha_{4}+\alpha_{5}+\alpha_{2},\quad    \alpha_{3}+\alpha_{4}+\alpha_{5}+\alpha_{6},\quad    \alpha_{4}+\alpha_{5}+\alpha_{6}+\alpha_{2},\quad    \alpha_{4}+\alpha_{5}+\alpha_{6}+\alpha_{7},\\
&  \alpha_{3}+2\alpha_{4}+\alpha_{5}+\alpha_{2},\quad    \alpha_{3}+\alpha_{4}+\alpha_{5}+\alpha_{6}+\alpha_{2},\quad    \alpha_{3}+\alpha_{4}+\alpha_{5}+\alpha_{6}+\alpha_{7},\\
& \alpha_{4}+\alpha_{5}+\alpha_{6}+\alpha_{7}+\alpha_{2},\quad    \alpha_{3}+2\alpha_{4}+\alpha_{5}+\alpha_{6}+\alpha_{2},\quad    \alpha_{3}+\alpha_{4}+\alpha_{5}+\alpha_{6}+\alpha_{7}+\alpha_{2},\\
&   \alpha_{3}+2\alpha_{4}+2\alpha_{5}+\alpha_{6}+\alpha_{2},\quad    \alpha_{3}+2\alpha_{4}+\alpha_{5}+\alpha_{6}+\alpha_{7}+\alpha_{2},\\
&   \alpha_{3}+2\alpha_{4}+2\alpha_{5}+\alpha_{6}+\alpha_{7}+\alpha_{2},
\quad    \alpha_{3}+2\alpha_{4}+2\alpha_{5}+2\alpha_{6}+\alpha_{7}+\alpha_{2}     \}
\eale
$

$\bale
\Delta_2^+=\{ &
\alpha_{1},\quad    \alpha_{4},\quad    \alpha_{6},\quad    \alpha_{4}+\alpha_{5},\quad    \alpha_{4}+\alpha_{2},\quad    \alpha_{5}+\alpha_{6},\quad    \alpha_{6}+\alpha_{7},\\
&     \alpha_{4}+\alpha_{5}+\alpha_{2},\quad    \alpha_{4}+\alpha_{5}+\alpha_{6},\quad    \alpha_{5}+\alpha_{6}+\alpha_{7},\\
&    \alpha_{4}+\alpha_{5}+\alpha_{6}+\alpha_{2},\quad    \alpha_{4}+\alpha_{5}+\alpha_{6}+\alpha_{7},\quad    \alpha_{4}+\alpha_{5}+\alpha_{6}+\alpha_{7}+\alpha_{2}    \}
\eale
$

$
\Delta_3^+=\
\alpha_{1},\quad    \alpha_{3},\quad    \alpha_{6},\quad    \alpha_{1}+\alpha_{3},\quad    \alpha_{5}+\alpha_{6},\quad    \alpha_{6}+\alpha_{7},\quad    \alpha_{5}+\alpha_{6}+\alpha_{7}     \}\\
$

$\bale
\Delta_4^+=\{&
\alpha_{1},\quad    \alpha_{3},\quad    \alpha_{4},\quad    \alpha_{1}+\alpha_{3},\quad    \alpha_{3}+\alpha_{4},\quad    \alpha_{4}+\alpha_{5},\quad    \alpha_{4}+\alpha_{2},\\
&  \alpha_{1}+\alpha_{3}+\alpha_{4},\quad    \alpha_{3}+\alpha_{4}+\alpha_{2},\quad    \alpha_{3}+\alpha_{4}+\alpha_{5},\quad    \alpha_{4}+\alpha_{5}+\alpha_{2},\\
&  \alpha_{1}+\alpha_{3}+\alpha_{4}+\alpha_{2},\quad    \alpha_{1}+\alpha_{3}+\alpha_{4}+\alpha_{5},\quad    \alpha_{3}+\alpha_{4}+\alpha_{5}+\alpha_{2},\quad    \alpha_{3}+2\alpha_{4}+\alpha_{5}+\alpha_{2},\\
&  \alpha_{1}+\alpha_{3}+\alpha_{4}+\alpha_{5}+\alpha_{2},\quad    \alpha_{1}+\alpha_{3}+2\alpha_{4}+\alpha_{5}+\alpha_{2},\quad    \alpha_{1}+2\alpha_{3}+2\alpha_{4}+\alpha_{5}+\alpha_{2}    \}.
\eale
$
\par
We get $s_k$ as follows,
$$
s_1=31, s_2=17, s_3=11, s_4=22.
$$
Then $s=31$.
\par

(15)EVII,$\mfg^\mbbc=\mfe_l, l=7,r=3$. \par
we have
$$
\al_1'=\la_1, \al_6'=\la_2,\al_7'=\la_3,
\al_2'=\al_3'=\al_4'=\al_5'=0.
$$\par
$\al'=0$ implies $m_2(\al)=m_3(\al)=m_4(\al)=m_5(\al)=0$.\par
We list $\Delta_k^+$ as follows,\\
$\bale
\quad \Delta_1^+=\{ &
\alpha_{6},\quad    \alpha_{7},\quad    \alpha_{5}+\alpha_{6},\quad    \alpha_{6}+\alpha_{7},\quad    \alpha_{4}+\alpha_{5}+\alpha_{6},\quad    \alpha_{5}+\alpha_{6}+\alpha_{7},\\
&   \alpha_{3}+\alpha_{4}+\alpha_{5}+\alpha_{6},\quad    \alpha_{4}+\alpha_{5}+\alpha_{6}+\alpha_{2},\quad    \alpha_{4}+\alpha_{5}+\alpha_{6}+\alpha_{7},\\
&   \alpha_{3}+\alpha_{4}+\alpha_{5}+\alpha_{6}+\alpha_{2},\quad    \alpha_{3}+\alpha_{4}+\alpha_{5}+\alpha_{6}+\alpha_{7},\quad    \alpha_{4}+\alpha_{5}+\alpha_{6}+\alpha_{7}+\alpha_{2},\\
&  \alpha_{3}+2\alpha_{4}+\alpha_{5}+\alpha_{6}+\alpha_{2},\quad    \alpha_{3}+\alpha_{4}+\alpha_{5}+\alpha_{6}+\alpha_{7}+\alpha_{2},\quad    \alpha_{3}+2\alpha_{4}+2\alpha_{5}+\alpha_{6}+\alpha_{2},\\
& \alpha_{3}+2\alpha_{4}+\alpha_{5}+\alpha_{6}+\alpha_{7}+\alpha_{2},\quad    \alpha_{3}+2\alpha_{4}+2\alpha_{5}+\alpha_{6}+\alpha_{7}+\alpha_{2},\\
& \alpha_{3}+2\alpha_{4}+2\alpha_{5}+2\alpha_{6}+\alpha_{7}+\alpha_{2}   \}
\eale
$

$\bale
\Delta_2^+=\{ &
\alpha_{1},\quad    \alpha_{7},\quad    \alpha_{1}+\alpha_{3},\quad    \alpha_{1}+\alpha_{3}+\alpha_{4},\\
&    \alpha_{1}+\alpha_{3}+\alpha_{4}+\alpha_{2},\quad    \alpha_{1}+\alpha_{3}+\alpha_{4}+\alpha_{5},\\
 & \alpha_{1}+\alpha_{3}+\alpha_{4}+\alpha_{5}+\alpha_{2},\quad    \alpha_{1}+\alpha_{3}+2\alpha_{4}+\alpha_{5}+\alpha_{2},\quad    \alpha_{1}+2\alpha_{3}+2\alpha_{4}+\alpha_{5}+\alpha_{2}     \}
\eale
$

$\bale
\Delta_3^+=\{&
\alpha_{1},\quad    \alpha_{6},\quad    \alpha_{1}+\alpha_{3},\quad    \alpha_{5}+\alpha_{6},\quad    \alpha_{1}+\alpha_{3}+\alpha_{4},\quad    \alpha_{4}+\alpha_{5}+\alpha_{6},\\
&   \alpha_{1}+\alpha_{3}+\alpha_{4}+\alpha_{2},\quad    \alpha_{1}+\alpha_{3}+\alpha_{4}+\alpha_{5},\quad    \alpha_{3}+\alpha_{4}+\alpha_{5}+\alpha_{6},\\
&   \alpha_{4}+\alpha_{5}+\alpha_{6}+\alpha_{2},\quad    \alpha_{1}+\alpha_{3}+\alpha_{4}+\alpha_{5}+\alpha_{2},\quad    \alpha_{1}+\alpha_{3}+\alpha_{4}+\alpha_{5}+\alpha_{6},\\
&  \alpha_{3}+\alpha_{4}+\alpha_{5}+\alpha_{6}+\alpha_{2},\quad    \alpha_{1}+\alpha_{3}+2\alpha_{4}+\alpha_{5}+\alpha_{2},\quad    \alpha_{1}+\alpha_{3}+\alpha_{4}+\alpha_{5}+\alpha_{6}+\alpha_{2},\\
& \alpha_{3}+2\alpha_{4}+\alpha_{5}+\alpha_{6}+\alpha_{2},\quad    \alpha_{1}+2\alpha_{3}+2\alpha_{4}+\alpha_{5}+\alpha_{2},\quad    \alpha_{1}+\alpha_{3}+2\alpha_{4}+\alpha_{5}+\alpha_{6}+\alpha_{2},\\
&   \alpha_{3}+2\alpha_{4}+2\alpha_{5}+\alpha_{6}+\alpha_{2}, \hq  \alpha_{1}+2\alpha_{3}+2\alpha_{4}+\alpha_{5}+\alpha_{6}+\alpha_{2}, \hq \alpha_{1}+\alpha_{3}+2\alpha_{4}+2\alpha_{5}+\alpha_{6}+\alpha_{2},\\
&  \alpha_{1}+2\alpha_{3}+2\alpha_{4}+2\alpha_{5}+\alpha_{6}+\alpha_{2},   \q \alpha_{1}+2\alpha_{3}+3\alpha_{4}+2\alpha_{5}+\alpha_{6}+\alpha_{2}, \\
& \alpha_{1}+2\alpha_{3}+3\alpha_{4}+2\alpha_{5}+\alpha_{6}+2\alpha_{2}   \}.
\eale
$\par
We get $s_k$ as follows,
$$
s_1=21, s_2=12, s_3=27.
$$
Then $s=27$.
\par

\vskip 1cm
We draw the Dynkin diagram  of  $\mfe_8$ as follows. \par
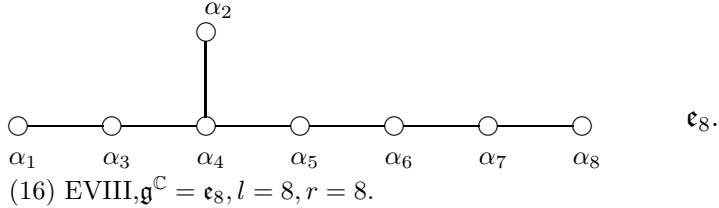
\begin{picture}(9, 3)
\put(0.125, 0.5){\circle{0.25}}
\put(0.25, 0.5){\line(1, 0){1}}
\put(1.375, 0.5){\circle{0.25}}
\put(1.5, 0.5){\line(1, 0){1.}}
\put(2.625, 0.5){\circle{0.25}}
\put(2.75, 0.5){\line(1, 0){1.}}
\put(3.875, 0.5){\circle{0.25}}
\put(4.0, 0.5){\line(1, 0){1.}}
\put(5.125, 0.5){\circle{0.25}}
\put(5.25, 0.5){\line(1, 0){1.}}
\put(6.375, 0.5){\circle{0.25}}
\put(6.5, 0.5){\line(1, 0){1.}}
\put(7.625, 0.5){\circle{0.25}}
\put(0, 0){$\alpha_1$}
\put(1.25, 0){$\alpha_3$}
\put(2.5, 0){$\alpha_4$}
\put(3.75, 0){$\alpha_5$}
\put(5, 0){$\alpha_6$}
\put(6.25, 0){$\alpha_7$}
\put(7.5, 0){$\alpha_8$}
\put(2.5, 0.625){ \line(0, 1){1} }
\put(2.625, 1.75){\circle{0.25}}
\put(2.6, 2){$\alpha_2$}
\put(9, 0.5){\Large{$\mathfrak{e}_8.$}}
\end{picture}
\par

(16) EVIII,$\mfg^\mbbc=\mfe_8, l=8,r=8$. \par
we have
$$
\al_i'=\la_i, 1\le i\le 8.
$$\par
$\al'=0$ implies $\al=0$. \par
We list $\Delta_k^+$ as follows,\\
$\bale \q
\Delta_1^+=\{&
\alpha_{3},\quad    \alpha_{4},\quad    \alpha_{5},\quad    \alpha_{6},\quad    \alpha_{7},\quad    \alpha_{8},\quad    \alpha_{2},\\
&   \alpha_{3}+\alpha_{4},\quad    \alpha_{4}+\alpha_{2},\quad    \alpha_{4}+\alpha_{5},\quad    \alpha_{5}+\alpha_{6},\quad    \alpha_{6}+\alpha_{7},\quad    \alpha_{7}+\alpha_{8},\\
& \alpha_{3}+\alpha_{4}+\alpha_{2},\quad    \alpha_{3}+\alpha_{4}+\alpha_{5},\quad    \alpha_{4}+\alpha_{5}+\alpha_{2},\\
&   \alpha_{4}+\alpha_{5}+\alpha_{6},\quad    \alpha_{5}+\alpha_{6}+\alpha_{7},\quad    \alpha_{6}+\alpha_{7}+\alpha_{8},\\
 & \alpha_{3}+\alpha_{4}+\alpha_{5}+\alpha_{2},\quad    \alpha_{3}+\alpha_{4}+\alpha_{5}+\alpha_{6},\quad    \alpha_{4}+\alpha_{5}+\alpha_{6}+\alpha_{2},\\
 &    \alpha_{4}+\alpha_{5}+\alpha_{6}+\alpha_{7},\quad    \alpha_{5}+\alpha_{6}+\alpha_{7}+\alpha_{8},\quad    \alpha_{3}+2\alpha_{4}+\alpha_{5}+\alpha_{2},\\
 & \alpha_{3}+\alpha_{4}+\alpha_{5}+\alpha_{6}+\alpha_{2},\quad    \alpha_{3}+\alpha_{4}+\alpha_{5}+\alpha_{6}+\alpha_{7},\quad    \alpha_{4}+\alpha_{5}+\alpha_{6}+\alpha_{7}+\alpha_{2},\\
  &  \alpha_{4}+\alpha_{5}+\alpha_{6}+\alpha_{7}+\alpha_{8},\quad    \alpha_{3}+2\alpha_{4}+\alpha_{5}+\alpha_{6}+\alpha_{2},\quad    \alpha_{3}+\alpha_{4}+\alpha_{5}+\alpha_{6}+\alpha_{7}+\alpha_{2},\\
  & \alpha_{3}+\alpha_{4}+\alpha_{5}+\alpha_{6}+\alpha_{7}+\alpha_{8},\quad    \alpha_{4}+\alpha_{5}+\alpha_{6}+\alpha_{7}+\alpha_{8}+\alpha_{2},\\
  &  \alpha_{3}+2\alpha_{4}+2\alpha_{5}+\alpha_{6}+\alpha_{2},
    \alpha_{3}+2\alpha_{4}+\alpha_{5}+\alpha_{6}+\alpha_{7}+\alpha_{2},\\
  &   \alpha_{3}+\alpha_{4}+\alpha_{5}+\alpha_{6}+\alpha_{7}+\alpha_{8}+\alpha_{2},\quad    \alpha_{3}+2\alpha_{4}+2\alpha_{5}+\alpha_{6}+\alpha_{7}+\alpha_{2},\\
  & \alpha_{3}+2\alpha_{4}+\alpha_{5}+\alpha_{6}+\alpha_{7}+\alpha_{8}+\alpha_{2},\quad    \alpha_{3}+2\alpha_{4}+2\alpha_{5}+2\alpha_{6}+\alpha_{7}+\alpha_{2},\\
  &   \alpha_{3}+2\alpha_{4}+2\alpha_{5}+\alpha_{6}+\alpha_{7}+\alpha_{8}+\alpha_{2},\quad    \alpha_{3}+2\alpha_{4}+2\alpha_{5}+2\alpha_{6}+\alpha_{7}+\alpha_{8}+\alpha_{2},\\
  &  \alpha_{3}+2\alpha_{4}+2\alpha_{5}+2\alpha_{6}+2\alpha_{7}+\alpha_{8}+\alpha_{2}     \}
  \eale
$

$\bale
\Delta_2^+=\{&
\alpha_{1},\quad    \alpha_{3},\quad    \alpha_{4},\quad    \alpha_{5},\quad    \alpha_{6},\quad    \alpha_{7},\quad    \alpha_{8},\\
&    \alpha_{1}+\alpha_{3},\quad    \alpha_{3}+\alpha_{4},\quad    \alpha_{4}+\alpha_{5},\quad    \alpha_{5}+\alpha_{6},\quad    \alpha_{6}+\alpha_{7},\quad    \alpha_{7}+\alpha_{8},\\
&  \alpha_{1}+\alpha_{3}+\alpha_{4},\quad    \alpha_{3}+\alpha_{4}+\alpha_{5},\quad    \alpha_{4}+\alpha_{5}+\alpha_{6},\\
&    \alpha_{5}+\alpha_{6}+\alpha_{7},\quad    \alpha_{6}+\alpha_{7}+\alpha_{8},\quad    \alpha_{1}+\alpha_{3}+\alpha_{4}+\alpha_{5},\\
&   \alpha_{3}+\alpha_{4}+\alpha_{5}+\alpha_{6},\quad    \alpha_{4}+\alpha_{5}+\alpha_{6}+\alpha_{7},\quad    \alpha_{5}+\alpha_{6}+\alpha_{7}+\alpha_{8},\\
& \alpha_{1}+\alpha_{3}+\alpha_{4}+\alpha_{5}+\alpha_{6},\quad    \alpha_{3}+\alpha_{4}+\alpha_{5}+\alpha_{6}+\alpha_{7},\quad    \alpha_{4}+\alpha_{5}+\alpha_{6}+\alpha_{7}+\alpha_{8},\\
& \alpha_{1}+\alpha_{3}+\alpha_{4}+\alpha_{5}+\alpha_{6}+\alpha_{7},\quad    \alpha_{3}+\alpha_{4}+\alpha_{5}+\alpha_{6}+\alpha_{7}+\alpha_{8},\\
&  \alpha_{1}+\alpha_{3}+\alpha_{4}+\alpha_{5}+\alpha_{6}+\alpha_{7}+\alpha_{8}     \}
\eale
$

$\bale
\Delta_3^+=\{ &
\alpha_{1},\quad    \alpha_{4},\quad    \alpha_{5},\quad    \alpha_{6},\quad    \alpha_{7},\quad    \alpha_{8},\quad    \alpha_{2},\\
&   \alpha_{4}+\alpha_{2},\quad    \alpha_{4}+\alpha_{5},\quad    \alpha_{5}+\alpha_{6},\quad    \alpha_{6}+\alpha_{7},\quad    \alpha_{7}+\alpha_{8},\\
&   \alpha_{4}+\alpha_{5}+\alpha_{2},\quad    \alpha_{4}+\alpha_{5}+\alpha_{6},\quad    \alpha_{5}+\alpha_{6}+\alpha_{7},\quad    \alpha_{6}+\alpha_{7}+\alpha_{8},\\
&   \alpha_{4}+\alpha_{5}+\alpha_{6}+\alpha_{2},\quad    \alpha_{4}+\alpha_{5}+\alpha_{6}+\alpha_{7},\quad    \alpha_{5}+\alpha_{6}+\alpha_{7}+\alpha_{8},\\
&  \alpha_{4}+\alpha_{5}+\alpha_{6}+\alpha_{7}+\alpha_{2},\quad    \alpha_{4}+\alpha_{5}+\alpha_{6}+\alpha_{7}+\alpha_{8},\quad    \alpha_{4}+\alpha_{5}+\alpha_{6}+\alpha_{7}+\alpha_{8}+\alpha_{2}     \}
\eale
$

$\bale
\Delta_4^+=\{ &
\alpha_{1},\quad    \alpha_{3},\quad    \alpha_{5},\quad    \alpha_{6},\quad    \alpha_{7},\quad    \alpha_{8},\quad    \alpha_{2},\\
&quad    \alpha_{1}+\alpha_{3},\quad    \alpha_{5}+\alpha_{6},\quad    \alpha_{6}+\alpha_{7},\quad    \alpha_{7}+\alpha_{8},\\
&  \alpha_{5}+\alpha_{6}+\alpha_{7},\quad    \alpha_{6}+\alpha_{7}+\alpha_{8},\quad    \alpha_{5}+\alpha_{6}+\alpha_{7}+\alpha_{8}     \}
\eale
$

$\bale
\Delta_5^+=\{ &
\alpha_{1},\quad    \alpha_{3},\quad    \alpha_{4},\quad    \alpha_{6},\quad    \alpha_{7},\quad    \alpha_{8},\quad    \alpha_{2},\\
& \alpha_{1}+\alpha_{3},\quad    \alpha_{3}+\alpha_{4},\quad    \alpha_{4}+\alpha_{2},\quad    \alpha_{6}+\alpha_{7},\quad    \alpha_{7}+\alpha_{8},\\
&  \alpha_{1}+\alpha_{3}+\alpha_{4},\quad    \alpha_{3}+\alpha_{4}+\alpha_{2},\quad    \alpha_{6}+\alpha_{7}+\alpha_{8},\\
&  \alpha_{1}+\alpha_{3}+\alpha_{4}+\alpha_{2}     \}
\eale
$

$\bale
\Delta_6^+=\{&
\alpha_{1},\quad    \alpha_{3},\quad    \alpha_{4},\quad    \alpha_{5},\quad    \alpha_{7},\quad    \alpha_{8},\quad    \alpha_{2},\\
 &  \alpha_{1}+\alpha_{3},\quad    \alpha_{3}+\alpha_{4},\quad    \alpha_{4}+\alpha_{2},\quad    \alpha_{4}+\alpha_{5},\quad    \alpha_{7}+\alpha_{8},\\
  &  \alpha_{1}+\alpha_{3}+\alpha_{4},\quad    \alpha_{3}+\alpha_{4}+\alpha_{2},\quad    \alpha_{3}+\alpha_{4}+\alpha_{5},\quad    \alpha_{4}+\alpha_{5}+\alpha_{2},\\
  & \alpha_{1}+\alpha_{3}+\alpha_{4}+\alpha_{5},\quad    \alpha_{1}+\alpha_{3}+\alpha_{4}+\alpha_{2},\quad    \alpha_{3}+\alpha_{4}+\alpha_{5}+\alpha_{2},\\
  &   \alpha_{1}+\alpha_{3}+\alpha_{4}+\alpha_{5}+\alpha_{2},\quad    \alpha_{3}+2\alpha_{4}+\alpha_{5}+\alpha_{2},\quad    \alpha_{1}+\alpha_{3}+2\alpha_{4}+\alpha_{5}+\alpha_{2},\\
  &   \alpha_{1}+2\alpha_{3}+2\alpha_{4}+\alpha_{5}+\alpha_{2}    \}
  \eale
$

$\bale
\Delta_7^+=\{ &
\alpha_{1},\quad    \alpha_{3},\quad    \alpha_{4},\quad    \alpha_{5},\quad    \alpha_{6},\quad    \alpha_{8},\quad    \alpha_{2},\\
&  \alpha_{1}+\alpha_{3},\quad    \alpha_{3}+\alpha_{4},\quad    \alpha_{4}+\alpha_{2},\quad    \alpha_{4}+\alpha_{5},\quad    \alpha_{5}+\alpha_{6},\\
& \alpha_{1}+\alpha_{3}+\alpha_{4},\quad    \alpha_{3}+\alpha_{4}+\alpha_{2},\quad    \alpha_{3}+\alpha_{4}+\alpha_{5},\\
 &    \alpha_{4}+\alpha_{5}+\alpha_{2},\quad    \alpha_{4}+\alpha_{5}+\alpha_{6},\\
 &    \alpha_{1}+\alpha_{3}+\alpha_{4}+\alpha_{5},\quad    \alpha_{1}+\alpha_{3}+\alpha_{4}+\alpha_{2},\quad    \alpha_{3}+\alpha_{4}+\alpha_{5}+\alpha_{2},\\
 &    \alpha_{3}+\alpha_{4}+\alpha_{5}+\alpha_{6},\quad    \alpha_{4}+\alpha_{5}+\alpha_{6}+\alpha_{2},\\
 &   \alpha_{1}+\alpha_{3}+\alpha_{4}+\alpha_{5}+\alpha_{2},\quad    \alpha_{3}+2\alpha_{4}+\alpha_{5}+\alpha_{2},\quad    \alpha_{1}+\alpha_{3}+\alpha_{4}+\alpha_{5}+\alpha_{6},\\
 &    \alpha_{3}+\alpha_{4}+\alpha_{5}+\alpha_{6}+\alpha_{2},\quad    \alpha_{1}+\alpha_{3}+2\alpha_{4}+\alpha_{5}+\alpha_{2},\quad    \alpha_{1}+\alpha_{3}+\alpha_{4}+\alpha_{5}+\alpha_{6}+\alpha_{2},\\
 &    \alpha_{3}+2\alpha_{4}+\alpha_{5}+\alpha_{6}+\alpha_{2},\quad    \alpha_{1}+2\alpha_{3}+2\alpha_{4}+\alpha_{5}+\alpha_{2},\quad    \alpha_{1}+\alpha_{3}+2\alpha_{4}+\alpha_{5}+\alpha_{6}+\alpha_{2},\\
 &   \alpha_{3}+2\alpha_{4}+2\alpha_{5}+\alpha_{6}+\alpha_{2},\quad    \alpha_{1}+2\alpha_{3}+2\alpha_{4}+\alpha_{5}+\alpha_{6}+\alpha_{2},\\
 &    \alpha_{1}+\alpha_{3}+2\alpha_{4}+2\alpha_{5}+\alpha_{6}+\alpha_{2},
    \alpha_{1}+2\alpha_{3}+2\alpha_{4}+2\alpha_{5}+\alpha_{6}+\alpha_{2},\\
 &    \alpha_{1}+2\alpha_{3}+3\alpha_{4}+2\alpha_{5}+\alpha_{6}+\alpha_{2},\quad    \alpha_{1}+2\alpha_{3}+3\alpha_{4}+2\alpha_{5}+\alpha_{6}+2\alpha_{2}     \}
 \eale
$

$\bale
\Delta_8^+=\{&
\alpha_{1},\quad    \alpha_{3},\quad    \alpha_{4},\quad    \alpha_{5},\quad    \alpha_{6},\quad    \alpha_{7},\quad    \alpha_{2},\\
 & \alpha_{1}+\alpha_{3},\quad    \alpha_{3}+\alpha_{4},\quad    \alpha_{4}+\alpha_{2},\quad    \alpha_{4}+\alpha_{5},\quad    \alpha_{5}+\alpha_{6},  \alpha_{6}+\alpha_{7},\\
 &    \alpha_{1}+\alpha_{3}+\alpha_{4},\quad    \alpha_{3}+\alpha_{4}+\alpha_{2},\quad    \alpha_{3}+\alpha_{4}+\alpha_{5},\\
 &   \alpha_{4}+\alpha_{5}+\alpha_{2},\quad    \alpha_{4}+\alpha_{5}+\alpha_{6},\quad    \alpha_{5}+\alpha_{6}+\alpha_{7},\\
 &   \alpha_{1}+\alpha_{3}+\alpha_{4}+\alpha_{5},\quad    \alpha_{1}+\alpha_{3}+\alpha_{4}+\alpha_{2},\quad    \alpha_{3}+\alpha_{4}+\alpha_{5}+\alpha_{2},\\
 &    \alpha_{3}+\alpha_{4}+\alpha_{5}+\alpha_{6},\quad    \alpha_{4}+\alpha_{5}+\alpha_{6}+\alpha_{2},\quad    \alpha_{4}+\alpha_{5}+\alpha_{6}+\alpha_{7},\\
 &    \alpha_{1}+\alpha_{3}+\alpha_{4}+\alpha_{5}+\alpha_{2},\quad    \alpha_{3}+2\alpha_{4}+\alpha_{5}+\alpha_{2},\quad    \alpha_{1}+\alpha_{3}+\alpha_{4}+\alpha_{5}+\alpha_{6},\\
 &   \alpha_{3}+\alpha_{4}+\alpha_{5}+\alpha_{6}+\alpha_{2},\quad    \alpha_{3}+\alpha_{4}+\alpha_{5}+\alpha_{6}+\alpha_{7},\quad    \alpha_{4}+\alpha_{5}+\alpha_{6}+\alpha_{7}+\alpha_{2},\\
 &    \alpha_{1}+\alpha_{3}+2\alpha_{4}+\alpha_{5}+\alpha_{2},\quad    \alpha_{1}+\alpha_{3}+\alpha_{4}+\alpha_{5}+\alpha_{6}+\alpha_{2},\quad    \alpha_{3}+2\alpha_{4}+\alpha_{5}+\alpha_{6}+\alpha_{2},\\
 &    \alpha_{1}+\alpha_{3}+\alpha_{4}+\alpha_{5}+\alpha_{6}+\alpha_{7},\quad    \alpha_{3}+\alpha_{4}+\alpha_{5}+\alpha_{6}+\alpha_{7}+\alpha_{2},\quad    \alpha_{1}+2\alpha_{3}+2\alpha_{4}+\alpha_{5}+\alpha_{2},\\
 &    \alpha_{1}+\alpha_{3}+2\alpha_{4}+\alpha_{5}+\alpha_{6}+\alpha_{2},\quad    \alpha_{3}+2\alpha_{4}+2\alpha_{5}+\alpha_{6}+\alpha_{2},\\
 & \alpha_{1}+\alpha_{3}+\alpha_{4}+\alpha_{5}+\alpha_{6}+\alpha_{7}+\alpha_{2}, \q
     \alpha_{3}+2\alpha_{4}+\alpha_{5}+\alpha_{6}+\alpha_{7}+\alpha_{2},\\
 &    \alpha_{1}+2\alpha_{3}+2\alpha_{4}+\alpha_{5}+\alpha_{6}+\alpha_{2},\quad    \alpha_{1}+\alpha_{3}+2\alpha_{4}+2\alpha_{5}+\alpha_{6}+\alpha_{2},\\
 &    \alpha_{1}+\alpha_{3}+2\alpha_{4}+\alpha_{5}+\alpha_{6}+\alpha_{7}+\alpha_{2},\quad    \alpha_{3}+2\alpha_{4}+2\alpha_{5}+\alpha_{6}+\alpha_{7}+\alpha_{2},\\
 &  \alpha_{1}+2\alpha_{3}+2\alpha_{4}+2\alpha_{5}+\alpha_{6}+\alpha_{2},\q
     \alpha_{1}+2\alpha_{3}+2\alpha_{4}+\alpha_{5}+\alpha_{6}+\alpha_{7}+\alpha_{2},\\
&    \alpha_{1}+\alpha_{3}+2\alpha_{4}+2\alpha_{5}+\alpha_{6}+\alpha_{7}+\alpha_{2},\quad    \alpha_{3}+2\alpha_{4}+2\alpha_{5}+2\alpha_{6}+\alpha_{7}+\alpha_{2},\\
 &   \alpha_{1}+2\alpha_{3}+3\alpha_{4}+2\alpha_{5}+\alpha_{6}+\alpha_{2},\quad    \alpha_{1}+2\alpha_{3}+2\alpha_{4}+2\alpha_{5}+\alpha_{6}+\alpha_{7}+\alpha_{2},\\
 &    \alpha_{1}+\alpha_{3}+2\alpha_{4}+2\alpha_{5}+2\alpha_{6}+\alpha_{7}+\alpha_{2},\quad    \alpha_{1}+2\alpha_{3}+3\alpha_{4}+2\alpha_{5}+\alpha_{6}+2\alpha_{2},\\
 &    \alpha_{1}+2\alpha_{3}+3\alpha_{4}+2\alpha_{5}+\alpha_{6}+\alpha_{7}+\alpha_{2},\quad    \alpha_{1}+2\alpha_{3}+2\alpha_{4}+2\alpha_{5}+2\alpha_{6}+\alpha_{7}+\alpha_{2},\\
 &  \alpha_{1}+2\alpha_{3}+3\alpha_{4}+2\alpha_{5}+\alpha_{6}+\alpha_{7}+2\alpha_{2},\quad    \alpha_{1}+2\alpha_{3}+3\alpha_{4}+2\alpha_{5}+2\alpha_{6}+\alpha_{7}+\alpha_{2},\\
 &    \alpha_{1}+2\alpha_{3}+3\alpha_{4}+2\alpha_{5}+2\alpha_{6}+\alpha_{7}+2\alpha_{2},\quad    \alpha_{1}+2\alpha_{3}+3\alpha_{4}+3\alpha_{5}+2\alpha_{6}+\alpha_{7}+\alpha_{2},\\
 &   \alpha_{1}+2\alpha_{3}+3\alpha_{4}+3\alpha_{5}+2\alpha_{6}+\alpha_{7}+2\alpha_{2},\quad    \alpha_{1}+2\alpha_{3}+4\alpha_{4}+3\alpha_{5}+2\alpha_{6}+\alpha_{7}+2\alpha_{2},\\
 &   \alpha_{1}+3\alpha_{3}+4\alpha_{4}+3\alpha_{5}+2\alpha_{6}+\alpha_{7}+2\alpha_{2},\quad    2\alpha_{1}+3\alpha_{3}+4\alpha_{4}+3\alpha_{5}+2\alpha_{6}+\alpha_{7}+2\alpha_{2}     \}.
\eale
$
\par
We get $s_k$ as follows,
$$
s_1=50, s_2=36, s_3=30,s_4=22,s_5=24, s_6=31, s_7=45,s_8=71
$$
Then $s=71$.
\par

(17) EIX,$\mfg^\mbbc=\mfe_8, l=8,r=4$. \par
we have
$$
\al_1'=\la_1, \al_6'=\la_2,\al_7'=\la_3,\al_8'=\la_4,
\al_2'=\al_3'=\al_4'=\al_5'=0.
$$\par
$\al'=0$ implies $m_2(\al)=m_3(\al)=m_4(\al)=m_5(\al)=0$.\par
We list $\Delta_k^+$ as follows,\\
$\bale
\Delta_1^+=\{ &
\alpha_{6},\quad    \alpha_{7},\quad    \alpha_{8},\quad    \alpha_{5}+\alpha_{6},\quad    \alpha_{6}+\alpha_{7},\quad    \alpha_{7}+\alpha_{8},\\
&    \alpha_{4}+\alpha_{5}+\alpha_{6},\quad    \alpha_{5}+\alpha_{6}+\alpha_{7},\quad    \alpha_{6}+\alpha_{7}+\alpha_{8},\\
&    \alpha_{3}+\alpha_{4}+\alpha_{5}+\alpha_{6},\quad    \alpha_{4}+\alpha_{5}+\alpha_{6}+\alpha_{2},\quad    \alpha_{4}+\alpha_{5}+\alpha_{6}+\alpha_{7},\\
& \alpha_{5}+\alpha_{6}+\alpha_{7}+\alpha_{8},\quad    \alpha_{3}+\alpha_{4}+\alpha_{5}+\alpha_{6}+\alpha_{2},\quad    \alpha_{3}+\alpha_{4}+\alpha_{5}+\alpha_{6}+\alpha_{7},\\
&    \alpha_{4}+\alpha_{5}+\alpha_{6}+\alpha_{7}+\alpha_{2},\quad    \alpha_{4}+\alpha_{5}+\alpha_{6}+\alpha_{7}+\alpha_{8},\quad    \alpha_{3}+2\alpha_{4}+\alpha_{5}+\alpha_{6}+\alpha_{2},\\
&    \alpha_{3}+\alpha_{4}+\alpha_{5}+\alpha_{6}+\alpha_{7}+\alpha_{2},\quad    \alpha_{3}+\alpha_{4}+\alpha_{5}+\alpha_{6}+\alpha_{7}+\alpha_{8},\quad    \alpha_{4}+\alpha_{5}+\alpha_{6}+\alpha_{7}+\alpha_{8}+\alpha_{2},\\
&    \alpha_{3}+2\alpha_{4}+2\alpha_{5}+\alpha_{6}+\alpha_{2},\quad    \alpha_{3}+2\alpha_{4}+\alpha_{5}+\alpha_{6}+\alpha_{7}+\alpha_{2},\\
&  \alpha_{3}+\alpha_{4}+\alpha_{5}+\alpha_{6}+\alpha_{7}+\alpha_{8}+\alpha_{2},\q
  \alpha_{3}+2\alpha_{4}+2\alpha_{5}+\alpha_{6}+\alpha_{7}+\alpha_{2},\\
  &  \alpha_{3}+2\alpha_{4}+\alpha_{5}+\alpha_{6}+\alpha_{7}+\alpha_{8}+\alpha_{2},\quad    \alpha_{3}+2\alpha_{4}+2\alpha_{5}+2\alpha_{6}+\alpha_{7}+\alpha_{2},\\
  &  \alpha_{3}+2\alpha_{4}+2\alpha_{5}+\alpha_{6}+\alpha_{7}+\alpha_{8}+\alpha_{2},\quad    \alpha_{3}+2\alpha_{4}+2\alpha_{5}+2\alpha_{6}+\alpha_{7}+\alpha_{8}+\alpha_{2},\\
  &   \alpha_{3}+2\alpha_{4}+2\alpha_{5}+2\alpha_{6}+2\alpha_{7}+\alpha_{8}+\alpha_{2}   \}
\\
\Delta_2^+=\{&
\alpha_{1},\quad    \alpha_{7},\quad    \alpha_{8},\quad    \alpha_{1}+\alpha_{3},\quad    \alpha_{7}+\alpha_{8},\q
   \alpha_{1}+\alpha_{3}+\alpha_{4},\\
   &    \alpha_{1}+\alpha_{3}+\alpha_{4}+\alpha_{5},\quad    \alpha_{1}+\alpha_{3}+\alpha_{4}+\alpha_{2},\quad    \alpha_{1}+\alpha_{3}+\alpha_{4}+\alpha_{5}+\alpha_{2},\\
   &  \alpha_{1}+\alpha_{3}+2\alpha_{4}+\alpha_{5}+\alpha_{2},\quad    \alpha_{1}+2\alpha_{3}+2\alpha_{4}+\alpha_{5}+\alpha_{2}     \}
\\
\Delta_3^+=\{&
\alpha_{1},\quad    \alpha_{6},\quad    \alpha_{8},\quad    \alpha_{1}+\alpha_{3},\quad    \alpha_{5}+\alpha_{6},\quad    \alpha_{1}+\alpha_{3}+\alpha_{4},\\
&    \alpha_{4}+\alpha_{5}+\alpha_{6},\quad    \alpha_{1}+\alpha_{3}+\alpha_{4}+\alpha_{5},\quad    \alpha_{1}+\alpha_{3}+\alpha_{4}+\alpha_{2},\\
&    \alpha_{3}+\alpha_{4}+\alpha_{5}+\alpha_{6},\quad    \alpha_{4}+\alpha_{5}+\alpha_{6}+\alpha_{2},\quad    \alpha_{1}+\alpha_{3}+\alpha_{4}+\alpha_{5}+\alpha_{2},\\
&    \alpha_{1}+\alpha_{3}+\alpha_{4}+\alpha_{5}+\alpha_{6},\quad    \alpha_{3}+\alpha_{4}+\alpha_{5}+\alpha_{6}+\alpha_{2},\quad    \alpha_{1}+\alpha_{3}+2\alpha_{4}+\alpha_{5}+\alpha_{2},\\
&    \alpha_{1}+\alpha_{3}+\alpha_{4}+\alpha_{5}+\alpha_{6}+\alpha_{2},\quad    \alpha_{3}+2\alpha_{4}+\alpha_{5}+\alpha_{6}+\alpha_{2},\\
&    \alpha_{1}+2\alpha_{3}+2\alpha_{4}+\alpha_{5}+\alpha_{2},\quad    \alpha_{1}+\alpha_{3}+2\alpha_{4}+\alpha_{5}+\alpha_{6}+\alpha_{2},\\
&    \alpha_{3}+2\alpha_{4}+2\alpha_{5}+\alpha_{6}+\alpha_{2},\quad    \alpha_{1}+2\alpha_{3}+2\alpha_{4}+\alpha_{5}+\alpha_{6}+\alpha_{2},\\
&   \alpha_{1}+\alpha_{3}+2\alpha_{4}+2\alpha_{5}+\alpha_{6}+\alpha_{2},\quad    \alpha_{1}+2\alpha_{3}+2\alpha_{4}+2\alpha_{5}+\alpha_{6}+\alpha_{2},\\
&    \alpha_{1}+2\alpha_{3}+3\alpha_{4}+2\alpha_{5}+\alpha_{6}+\alpha_{2},\quad    \alpha_{1}+2\alpha_{3}+3\alpha_{4}+2\alpha_{5}+\alpha_{6}+2\alpha_{2}     \}
\\
\Delta_4^+=\{&
\alpha_{1},\quad    \alpha_{6},\quad    \alpha_{7},\quad    \alpha_{1}+\alpha_{3},\quad    \alpha_{5}+\alpha_{6},\quad    \alpha_{6}+\alpha_{7},\\
&   \alpha_{1}+\alpha_{3}+\alpha_{4},\quad    \alpha_{4}+\alpha_{5}+\alpha_{6},\quad    \alpha_{5}+\alpha_{6}+\alpha_{7},\\
&   \alpha_{1}+\alpha_{3}+\alpha_{4}+\alpha_{5},\quad    \alpha_{1}+\alpha_{3}+\alpha_{4}+\alpha_{2},\quad    \alpha_{3}+\alpha_{4}+\alpha_{5}+\alpha_{6},\\
&    \alpha_{4}+\alpha_{5}+\alpha_{6}+\alpha_{2},\quad    \alpha_{4}+\alpha_{5}+\alpha_{6}+\alpha_{7},\quad    \alpha_{1}+\alpha_{3}+\alpha_{4}+\alpha_{5}+\alpha_{2},\\
&    \alpha_{1}+\alpha_{3}+\alpha_{4}+\alpha_{5}+\alpha_{6},\quad    \alpha_{3}+\alpha_{4}+\alpha_{5}+\alpha_{6}+\alpha_{2},\quad    \alpha_{3}+\alpha_{4}+\alpha_{5}+\alpha_{6}+\alpha_{7},\\
&   \alpha_{4}+\alpha_{5}+\alpha_{6}+\alpha_{7}+\alpha_{2},\quad    \alpha_{1}+\alpha_{3}+2\alpha_{4}+\alpha_{5}+\alpha_{2},\\
&    \alpha_{1}+\alpha_{3}+\alpha_{4}+\alpha_{5}+\alpha_{6}+\alpha_{2},\quad    \alpha_{3}+2\alpha_{4}+\alpha_{5}+\alpha_{6}+\alpha_{2},\\
&   \alpha_{1}+\alpha_{3}+\alpha_{4}+\alpha_{5}+\alpha_{6}+\alpha_{7},\quad    \alpha_{3}+\alpha_{4}+\alpha_{5}+\alpha_{6}+\alpha_{7}+\alpha_{2},\\
&   \alpha_{1}+2\alpha_{3}+2\alpha_{4}+\alpha_{5}+\alpha_{2},\quad    \alpha_{1}+\alpha_{3}+2\alpha_{4}+\alpha_{5}+\alpha_{6}+\alpha_{2},\\
&    \alpha_{3}+2\alpha_{4}+2\alpha_{5}+\alpha_{6}+\alpha_{2},\quad    \alpha_{1}+\alpha_{3}+\alpha_{4}+\alpha_{5}+\alpha_{6}+\alpha_{7}+\alpha_{2},\\
&   \alpha_{3}+2\alpha_{4}+\alpha_{5}+\alpha_{6}+\alpha_{7}+\alpha_{2},\quad    \alpha_{1}+2\alpha_{3}+2\alpha_{4}+\alpha_{5}+\alpha_{6}+\alpha_{2},\\
&    \alpha_{1}+\alpha_{3}+2\alpha_{4}+2\alpha_{5}+\alpha_{6}+\alpha_{2},\quad    \alpha_{1}+\alpha_{3}+2\alpha_{4}+\alpha_{5}+\alpha_{6}+\alpha_{7}+\alpha_{2},\\
&    \alpha_{3}+2\alpha_{4}+2\alpha_{5}+\alpha_{6}+\alpha_{7}+\alpha_{2},\quad    \alpha_{1}+2\alpha_{3}+2\alpha_{4}+2\alpha_{5}+\alpha_{6}+\alpha_{2},\\
&    \alpha_{1}+2\alpha_{3}+2\alpha_{4}+\alpha_{5}+\alpha_{6}+\alpha_{7}+\alpha_{2},\quad    \alpha_{1}+\alpha_{3}+2\alpha_{4}+2\alpha_{5}+\alpha_{6}+\alpha_{7}+\alpha_{2},\\
&    \alpha_{3}+2\alpha_{4}+2\alpha_{5}+2\alpha_{6}+\alpha_{7}+\alpha_{2},\quad    \alpha_{1}+2\alpha_{3}+3\alpha_{4}+2\alpha_{5}+\alpha_{6}+\alpha_{2},\\
&    \alpha_{1}+2\alpha_{3}+2\alpha_{4}+2\alpha_{5}+\alpha_{6}+\alpha_{7}+\alpha_{2},\quad    \alpha_{1}+\alpha_{3}+2\alpha_{4}+2\alpha_{5}+2\alpha_{6}+\alpha_{7}+\alpha_{2},\\
&    \alpha_{1}+2\alpha_{3}+3\alpha_{4}+2\alpha_{5}+\alpha_{6}+2\alpha_{2},\quad    \alpha_{1}+2\alpha_{3}+3\alpha_{4}+2\alpha_{5}+\alpha_{6}+\alpha_{7}+\alpha_{2},\\
&    \alpha_{1}+2\alpha_{3}+2\alpha_{4}+2\alpha_{5}+2\alpha_{6}+\alpha_{7}+\alpha_{2},\quad    \alpha_{1}+2\alpha_{3}+3\alpha_{4}+2\alpha_{5}+\alpha_{6}+\alpha_{7}+2\alpha_{2},\\
&   \alpha_{1}+2\alpha_{3}+3\alpha_{4}+2\alpha_{5}+2\alpha_{6}+\alpha_{7}+\alpha_{2},\quad    \alpha_{1}+2\alpha_{3}+3\alpha_{4}+2\alpha_{5}+2\alpha_{6}+\alpha_{7}+2\alpha_{2},\\
&    \alpha_{1}+2\alpha_{3}+3\alpha_{4}+3\alpha_{5}+2\alpha_{6}+\alpha_{7}+\alpha_{2},\quad    \alpha_{1}+2\alpha_{3}+3\alpha_{4}+3\alpha_{5}+2\alpha_{6}+\alpha_{7}+2\alpha_{2},\\
&    \alpha_{1}+2\alpha_{3}+4\alpha_{4}+3\alpha_{5}+2\alpha_{6}+\alpha_{7}+2\alpha_{2},\quad    \alpha_{1}+3\alpha_{3}+4\alpha_{4}+3\alpha_{5}+2\alpha_{6}+\alpha_{7}+2\alpha_{2},\\
&   2\alpha_{1}+3\alpha_{3}+4\alpha_{4}+3\alpha_{5}+2\alpha_{6}+\alpha_{7}+2\alpha_{2}   \}.
\eale
$\

We get $s_k$ as follows,
$$
s_1=34, s_2=15, s_3=29,s_4=55.
$$
Then $s=55$.
\par

\vskip 1cm
The Dynkin diagram of $\mathfrak{f}_4$ is \\
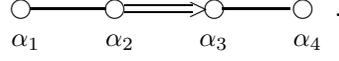
\centerline
{
\begin{picture}(9, 1)
\put(0.125, 0.5){\circle{0.25}}
\put(0.25, 0.5){\line(1, 0){1}}
\put(1.375, 0.5){\circle{0.25}}
\put(1.5, 0.55){\line(1, 0){0.9}}
\put(1.5, 0.45){\line(1, 0){0.9}}
\put(2.35, 0.40){$>$}
\put(2.70, 0.5){\circle{0.25}}
\put(2.8, 0.5){\line(1, 0){0.9}}
\put(3.875, 0.5){\circle{0.25}}
\put(0, 0){$\alpha_1$}
\put(1.25, 0){$\alpha_2$}
\put(2.5, 0){$\alpha_3$}
\put(3.75, 0){$\alpha_4$}
\put(4.3, 0.4){.}
\end{picture}
}
We list all positive roots
\begin{eqnarray*}
&& \alpha_1, \quad \alpha_2, \quad \alpha_3, \quad \alpha_4\\
&& \alpha_1+\alpha_2, \quad \alpha_2+\alpha_3, \quad \alpha_3+\alpha_4\\
&& 2\alpha_2+\alpha_3, \quad \alpha_1+\alpha_2+\alpha_3, \quad \alpha_2+\alpha_3+\alpha_4\\
&& \alpha_1+\alpha_2+2\alpha_3, \quad \alpha_1+\alpha_2+\alpha_3+\alpha_4, \quad \alpha_2+2\alpha_3+\alpha_4\\
&& \alpha_1+2\alpha_2+2\alpha_3, \quad \alpha_2+2\alpha_3+2\alpha_4, \quad \alpha_1+\alpha_2+2\alpha_3+\alpha_4\\
&& \alpha_1+\alpha_2+2\alpha_3+2\alpha_4, \quad \alpha_1+2\alpha_2+2\alpha_3+\alpha_4\\
&& \alpha_1+2\alpha_2+2\alpha_3+2\alpha_4, \quad \alpha_1+2\alpha_2+3\alpha_3+\alpha_4\\
&& \alpha_1+2\alpha_2+3\alpha_3+2\alpha_4\\
&& \alpha_1+2\alpha_2+4\alpha_3+2\alpha_4\\
&& \alpha_1+3\alpha_2+4\alpha_3+2\alpha_4\\
&& 2\alpha_1+3\alpha_2+4\alpha_3+2\alpha_4.
\end{eqnarray*}

(18) FI,$\mfg^\mbbc=\mff_4, l=4,r=4$. \par
we have
$$
\al_1'=\la_1, \al_2'=\la_2,\al_3'=\la_3,\al_4'=\la_4.
$$\par
$\al'=0$ implies $\al=0$.\par
We list $\Delta_k^+$ as follows,\\
$\bale
\q\Delta_1^+=\{&
\alpha_{2},\quad    \alpha_{3},\quad    \alpha_{4},\quad    \alpha_{2}+\alpha_{3},\quad    \alpha_{3}+\alpha_{4},\quad    \alpha_{2}+2\alpha_{3},\\
&    \alpha_{2}+\alpha_{3}+\alpha_{4},\quad    \alpha_{2}+2\alpha_{3}+\alpha_{4},\quad    \alpha_{2}+2\alpha_{3}+2\alpha_{4}     \}\\
\q\Delta_2^+=\{&
\alpha_{1},\quad    \alpha_{3},\quad    \alpha_{4},\quad    \alpha_{3}+\alpha_{4}   \}\\
\q\Delta_3^+=\{&
\alpha_{1},\quad    \alpha_{2},\quad    \alpha_{4},\quad    \alpha_{1}+\alpha_{2}    \}\\
\q\Delta_4^+=\{&
\alpha_{1},\quad    \alpha_{2},\quad    \alpha_{3},\quad    \alpha_{1}+\alpha_{2},\quad    \alpha_{2}+\alpha_{3},\\
&   \alpha_{2}+2\alpha_{3},\quad    \alpha_{1}+\alpha_{2}+\alpha_{3},\quad    \alpha_{1}+\alpha_{2}+2\alpha_{3},\quad    \alpha_{1}+2\alpha_{2}+2\alpha_{3}    \}.
\eale
$

\par
We get $s_k$ as follows,
$$
s_1=13, s_2=8, s_3=8,s_4=13.
$$
Then $s=13$.
\par

(19) FII,$\mfg^\mbbc=\mff_4, l=4,r=1$. \par
Then $s=r=1$.
\par

\vskip 1cm
The Dynkin diagram of $\mathfrak{g}_2$ is \\
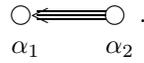
\centerline
{
\begin{picture}(9, 2)
\put(0.125, 0.5){\circle{0.25}}
\put(0.25, 0.40){$<$}
\put(0.35, 0.5){\line(1, 0){0.9}}
\put(0.35, 0.55){\line(1, 0){0.9}}
\put(0.35, 0.45){\line(1, 0){0.9}}
\put(1.375, 0.5){\circle{0.25}}
\put(0, 0){$\alpha_1$}
\put(1.25, 0){$\alpha_2$}
\put(1.7,0.4){.}
\end{picture}
}\\
We list all positive roots
\begin{eqnarray*}
&& \alpha_1, \quad \alpha_2\\
&& \alpha_1+\alpha_2, \quad  2\alpha_1+\alpha_2\\
&& 3\alpha_1+\alpha_2, \quad  3\alpha_1+2\alpha_2.
\end{eqnarray*}
\par

(20) G,$\mfg^\mbbc=\mfg_2, l=2,r=2$. \par
we have
$$
\al_1'=\la_1, \al_2'=\la_2.
$$\par
$\al'=0$ implies $\al=0$.\par
We list $\Delta_k^+$ as follows,\\
$
\Delta_1^+=\{
\alpha_{2}     \}\\
$\\
$
\Delta_2^+=\{
\alpha_{1}     \}.\\
$\par
We get $s_k$ as follows,
$$
s_1=3, s_2=3.
$$\par
Then $s=3$.
\par

\begin{remark}
For convenience, in the case of $EI$ we give another computations  by using
the embedding listed in \cite[p81]{Sam}.
\end{remark}
In $\mbbr^6$, the positive  roots of $\mfe_6$ are:
$$\vep_i-\vep_j, i<j;  \q \vep_i+\vep_j+\vep_k, i<j<k; \q \vep_1+\vep_2+\vep_3+\vep_4+\vep_5+\vep_6.$$\par
The simple roots are
$$
\al_i=\vep_i-\vep_{i+1}, 1\le i\le 5; \al_6=\vep_4+\vep_5+\vep_6.
$$
We have
\beqs
\vep_i-\vep_j\qq\ &=&\al_i+\al_{i-1}+\cdots +\al_{j-1}, \q 1\le i <j\le 6\\
\vep_1+\vep_2+\vep_3&=&\al_1+2\al_2+3\al_3+2\al_4+\al_5+\al_6\\
\vep_1+\vep_2+\vep_4&=&\al_1+2\al_2+2\al_3+2\al_4+\al_5+\al_6\\
\vep_1+\vep_2+\vep_5&=&\al_1+2\al_2+2\al_3+\al_4+\al_5+\al_6\\
\vep_1+\vep_2+\vep_6&=&\al_1+2\al_2+2\al_3+\al_4   +\al_6\\
\vep_1+\vep_3+\vep_4&=&\al_1+\al_2+2\al_3+2\al_4+\al_5+\al_6\\
\vep_1+\vep_3+\vep_5&=&\al_1+\al_2+2\al_3+\al_4+\al_5+\al_6\\
\vep_1+\vep_3+\vep_6&=&\al_1+\al_2+2\al_3+\al_4 +\al_6\\
\vep_1+\vep_4+\vep_5&=&\al_1+\al_2+\al_3+\al_4+\al_5+\al_6\\
\vep_1+\vep_4+\vep_6&=&\al_1+\al_2+\al_3+\al_4 +\al_6\\
\vep_1+\vep_5+\vep_6&=&\al_1+\al_2+\al_3  +\al_6\\
\vep_2+\vep_3+\vep_4&=&\al_2+2\al_3+2\al_4+\al_5+\al_6\\
\vep_2+\vep_3+\vep_5&=&\al_2+2\al_3+\al_4+\al_5+\al_6\\
\vep_2+\vep_3+\vep_6&=&\al_2+2\al_3+2\al_4 +\al_6\\
\vep_2+\vep_4+\vep_5&=&\al_2+\al_3+\al_4+\al_5+\al_6\\
\vep_2+\vep_4+\vep_6&=&\al_2+\al_3+\al_4 +\al_6\\
\vep_2+\vep_5+\vep_6&=&\al_2+\al_3 +\al_6\\
\vep_3+\vep_4+\vep_5&=&\al_3+\al_4+\al_5 +\al_6\\
\vep_3+\vep_4+\vep_6&=&\al_3+\al_4 +\al_6\\
\vep_3+\vep_5+\vep_6&=&\al_3 +\al_6\\
\vep_4+\vep_5+\vep_6&=&\al_6\\
\vep_1+\vep_2+\vep_3+\vep_4+\vep_5+\vep_6&=&\al_1+2\al_2+3\al_3+2\al_4+\al_5+2\al_6.
\eeqs
\par
To get $\Delta_k^+$, we have two cases.\par
For $\al=\vep_i-\vep_j, i<j, m_k'(\al)=0$, then  $ i<j\le k \mbox{ or }  j>i>k$, the number of roots  is
$$\binom{k}{2}+ \binom{6-k}{2}.$$\par
For $\al=\vep_i+\vep_j+\vep_p, i<j<p,  m_k'(\al)=0$, then we have \\
\beqs
&&k=1,2,3, \q  k<i<j<p, \mbox {the number of positive roots is }\binom{6-k}{3}.\\
&&k=4, \al_1+\al_5+\al_6, \al_2+\al_5+\al_6, \al_3+\al_5+\al_6, \al_4+\al_5+\al_6.\\
&&k=5, \al_1+\al_2+\al_6, \al_1+\al_3+\al_6,\al_1+\al_4+\al_6,\al_1+\al_5+\al_6,\al_2+\al_3+\al_6,\\
&& \q \q    \al_2+\al_4+\al_6,\al_2+\al_5+\al_6,\al_3+\al_4+\al_6,\al_3+\al_5+\al_6,\al_4+\al_5+\al_6.\\
&&k=6, \mbox{empty. }\\
\eeqs
\par
Then we get the values of $s_k$ are $26,17,13,17,26,21$ respectively.\par

\vskip 1cm
Now we get all s-values  for  irreducible Riemannian
symmetric spaces of compact type. 
In summary,  we have the  table 1.1 in the 
introduction of the present paper.
\par

\end{document}